\documentclass{article}

\usepackage[english]{babel}
\usepackage[T1]{fontenc}
\usepackage[utf8]{inputenc}
\usepackage{amsfonts}
\usepackage{amsmath}
\usepackage{dsfont}
\usepackage{hyperref}
\usepackage[thmmarks,amsmath]{ntheorem}
\usepackage{amssymb}
\usepackage{geometry}
\usepackage{esint}
\usepackage{enumitem}

\numberwithin{equation}{section}

{
\newtheorem{Def}{Définition}[section]
}

{
\newtheorem{Lem}[Def]{\textit{Lemma}}
}

{

}

\newtheorem{Prop}[Def]{\textit{Proposition}}

\newtheorem{Theo}[Def]{Theorem}

{
\theorembodyfont{}
\newtheorem{Rem}[Def]{Remark}
}

{
\theoremstyle{break}
\theoremsymbol{$\square$}
\theorembodyfont{}
\newtheorem*{Dem}{Proof}
}

\allowdisplaybreaks

\begin{document}

\title{Scattering of the 2D modified Zakharov-Kuznetsov equation}

\author{Philippe Anjolras}

\maketitle

\begin{abstract}
We study the modified Zakharov-Kuznetsov equation in dimension $2$ : 
\[ \partial_t u + \partial_x \left( \Delta u + u^3 \right) = 0 \]
where $u : (t, (x, y)) \in \mathbb{R} \times \mathbb{R}^2 \mapsto u(t, x, y) \in \mathbb{R}$ and $\Delta = \partial_x^2 + \partial_y^2$ is the full Laplacian. We prove that solutions for small and localized initial data scatter for large time. Our proof relies on the method of space-time resonances. 
\end{abstract}

\tableofcontents

\section{Introduction}

We study the modified Zakharov-Kuznetsov (mZK) equation in two spatial dimensions :  
\begin{equation} \left\{ \begin{array}{l} \partial_t u + \partial_x \left( \Delta u + \frac{u^3}{3} \right) = 0 \\
u(t=0) = u_0 \end{array} \right. \label{equ_ZK} \end{equation}
where $u : (t, x, y) \in \mathbb{R}^3 \mapsto u(t, x, y) \in \mathbb{R}$, and $\Delta = \partial_{xx} + \partial_{yy}$. The original equation derived by Zakharov and Kuznetsov (in dimension $d = 3$) has only a quadratic non-linearity $\frac{u^2}{2}$ \cite{ZKoriginal} and describes the propagation of ion-sound waves in a plasma ; it was derived rigorously from the Euler-Poisson system by Lannes, Linares and Saut in \cite{LLS-derivation}, while Chiron showed in \cite{Chironlimder} that in degenerate cases a cubic non-linearity could be obtained at the limit. Mathematically, this equation corresponds to a multi-dimensional extension of the (modified) Korteweg-de-Vries equation. 

\paragraph{Well-posedness of the modified Zakharov-Kuznetsov equation} The well-posedness of equation \eqref{equ_ZK} has been extensively studied in the past years. (mZK) enjoys the scaling invariance : 
\[ u_{\lambda}(t, x, y) = \lambda u(\lambda^3 t, \lambda x, \lambda y) \]
so that the critical scaling index is $s_c = 0$. However, Kinoshita proved in \cite{KinoshitaLWPmZK2D} that local well-posedness holds in $H^s(\mathbb{R}^2)$ for any $s \geq \frac{1}{4}$, while the data-to-solution map fails to be regular at lower regularity. This result extended previous works by Biagioni and Linares in $H^1$ \cite{BiagioniLinaresLWP}, Linares and Pastor in $H^{3/4+}$ \cite{LinaresPastorLWP}, Ribaud and Vento in $H^{1/4+}$ \cite{RibaudVentoLWP}. Local well-posedness results have also been proved for the Zakharov-Kuznetsov equation \cite{FaminskiiLWP, GrunrockHerrLWP, KinoshitaGWP, MolinetPilodLWP} or for the generalized equation, where the non-linearity is some $u^{1+k}$, $k \geq 3$ \cite{FarahLinaresgZK, GrunrockgZK, LinaresPastorgZK, RibaudVentoLWP} as well as in dimension $d = 3$ \cite{GrunrockLWP3d, GrunrockgZK, LinaresSautLWP, RibaudVentoLWP3d, HerrKinoshitaLWP}. 

Concerning the long-time behaviour, Linares and Pastor proved global well-posedness in $H^s$ for $s > \frac{53}{63}$ \cite{LinaresPastorgZK} and \cite{BFR-GWPmZK2D}, while scattering theory was only proved for the generalized Zakharov-Kuznetsov equation \cite{FarahLinaresgZK}, or in dimension $1$ for the mKdV equation \cite{HayashiNaumkinmKdV1, HayashiNaumkinmKdV2, HarropGriffLongtimemKdV, GermainPusateriRousset_mKDV}. On the other hand, the question of the (in)stability near a soliton has been investigated recently in \cite{Chenetalsoliton, Masmoudietalsoliton}. 

The purpose of this note is to prove a non-linear scattering result for small data. Our method is based on the space-time resonances analysis developped by Germain, Masmoudi and Shatah in \cite{Germain_3DSchrod, Germain_2DSchrod, GMSWaterwaves} for the quadratic Schrödinger equation and gravity water waves, and simultaneously by Gustafson, Nakanishi and Tsai in \cite{NakanishiResonances, NakanishiResonances2} for the Gross-Pitaevskii equation. The space-time resonance methode has been widely used in the recent years to prove global regularity of nonlinear PDEs, e.g. water waves \cite{GMSWaterwaves,IonescuPusateri_waterwaves,GMS_waterwavescap}, the Klein-Gordon equation \cite{IonescuPausader_KG}, the wave equation \cite{ShatahPusateri,AnjolrasABI}, the Euler-Maxwell system \cite{Germain_EulerMaxwell, GuoIonescuPausader_plasma,GuoIonescuPausader_EulerMaxwell,DengIonescuPausader_EulerMaxwell}, the Euler-Poisson system \cite{GuoPausader_EulerPoisson,GuoIonescuPausader_plasma}. It has also been used to get scattering \cite{IonescuPausader_KG,IonescuPusateri_waterwaves,GMS_waterwavescap,GermainPusateriRousset_mKDV,AnjolrasABI}. Similar results concerning global regularity and scattering have also been obtained in \cite{DelortAlazard_gravitywaterwaves2015,HunterIfrimTataru_ww,IfrimTataru_ww}. We explain the details of the method in section \ref{section-spacetimeresonance}. There were also recent advances for cubic dispersive equations by Ifrim and Tataru concerning the long-time behaviour of solutions in dimension $1$, see for instance \cite{IfrimTataru1Dconj}. 

\paragraph{} By the change of variables 
\[ (x_a, x_b) = \left( \frac{x + \sqrt{3} y}{2^{2/3}}, \frac{x - \sqrt{3} y}{2^{2/3}} \right), \quad v(t, x_a, x_b) = 2^{-1/2} u(t, x, y) \]
we may transform equation \eqref{equ_ZK} into 
\begin{equation} \partial_t v + \partial_{x_a}^3 v + \partial_{x_b}^3 v + (\partial_{x_a} + \partial_{x_b}) (v^3) = 0 \label{equ_ZK_coordab} \end{equation}

Denote by $L = \partial_{x_a}^3 + \partial_{x_b}^3$ the linear operator associated to \eqref{equ_ZK_coordab}. Recall that the profile function is defined as $f(t) = e^{tL} v(t)$, so that $f$ satisfies
\[ \partial_t f + e^{-tL} \left( \partial_{x_a} + \partial_{x_b} \right) \left( e^{tL} f \right)^3 = 0 \]
In what follows, we won't explicitely distinguish between $u$ and $v$, since they are related by a simple invertible change of variables. 

Our main theorem is : 

\begin{Theo} \label{main-theo} There exists $\varepsilon_0, C > 0$ such that, for any initial data $u_0$ satisfying 
\[ \Vert u_0 \Vert_{H^3} + \Vert (|x|+|y|) u_0 \Vert_{L^2} \leq \varepsilon_0 \]
then the global solution $u$ of \eqref{equ_ZK} with initial data $u_0$ satisfies 
\[ \Vert u(t) \Vert_{H^3} \leq C \varepsilon_0 \]
and the associated profile function $f(t)$ converges in $H^2$ as $t \to \infty$. 
\end{Theo}

\begin{Rem} \begin{itemize}
\item We will get an explicit bound $\Vert f(\infty) - f(0) \Vert_{H^2} \lesssim \varepsilon^3$. We focus on positive times, but it is clear that the results applies as well for $t \to -\infty$. 
\item One could expect to be able to lower the requirement of an $H^3$ bound to $H^1$, the energy space, with scattering in $L^2$. However, the $H^3$ requirement makes our dispersive estimate much easier : it requires only usual harmonic analysis and well-known kernel estimates. The $H^3$ requirement also simplifies the weighted estimate, although our proof could be applied with minor changes. Finally, one would need to deal with the short-time propagation of the weighted norm, which should be possible adapting proofs of local well-posedness, while we completely avoid this issue and focus on long-time behaviour. We therefore preferred to state and prove a slightly sub-optimal result to keep the attention on the long-time behaviour and the use of the non-resonant structure, in an approach that would be the least technical. 
\end{itemize}
\end{Rem}

\paragraph{Outline of the article and major steps of the proof} In section \ref{section-dispersive-estimate}, we recall results regarding the time decay of the kernel associated to the linear evolution, that will be used in section \ref{section-prelest} to prove decay in $L^{\infty}$ of the nonlinear solution. In section \ref{section-paraproducts}, we recall results regarding the continuity of singular paraproducts. We state the main a priori estimate (proposition \ref{prop-est-a-priori}) in section \ref{section-statement-apriori}, and prove our main theorem in section \ref{section-proof-mainth} assuming the a priori estimate holds. The rest of the paper is devoted to the proof of proposition \ref{prop-est-a-priori}. We prove a $H^3$ estimate in section \ref{section-H3est}, using standard energy identities. 

The last section \ref{section-weightedest} contains the proof of the weighted estimate, that is separated in three steps. We first perform a space-time resonances analysis in the spirit of Germain-Masmoudi-Shatah. In section \ref{section-genericweightedest}, we prove a generic estimate for weighted terms, and prove the weighted estimate in section \ref{section-proofweighted}. Although we only dispose of a singular non-resonant structure, the dispersive decay is enough to close the argument up to performing a time-dependent frequency separation removing the singular area. 

We now give a more detailed overview of this weighted estimate. Recall that the space-time resonance method relies on writing a Duhamel formula for the Fourier transform of the profile function $f(t) = e^{tL} u(t)$, so that 
\begin{equation} \widehat{f}(t, \xi) = \widehat{f}(0, \xi) + \int_0^t \int \int e^{i s \varphi(\xi, \eta, \sigma)} m(\xi, \eta, \sigma) \widehat{f}(s, \xi-\eta-\sigma) \widehat{f}(s, \eta) \widehat{f}(s, \sigma) ~ d\sigma d\eta ds \label{equDuhamelheuristics} \end{equation}
where $m$ is the symbol associated to the non-linearity (here $m(\xi, \eta, \sigma) = \xi_a + \xi_b$) while $\varphi$ is the interaction phase : 
\[ \varphi(\xi, \eta, \sigma) = \xi_a^3 + \xi_b^3 - (\xi_a - \eta_a - \sigma_a)^3 - (\xi_b - \eta_b - \sigma_b)^3 - \eta_a^3 - \eta_b^3 - \sigma_a^3 - \sigma_b^3 \]
From the stationary phase heuristics and the fact that, in the Duhamel formula \eqref{equDuhamelheuristics}, one has an integration in $s, \eta, \sigma$, the space-time resonance method relies on the heuristic that the dynamic should be driven by the behaviour near the resonant set
\[ \{ \varphi = 0 \} \cap \{ \nabla_{\eta} \varphi = \nabla_{\sigma} \varphi = 0 \} \]
In our case, since we will apply it to estimate a weighted norm on $f$, this translates into an estimate on $\nabla_{\xi} \widehat{f}(t)$. Therefore, distributing the derivative in \eqref{equDuhamelheuristics}, one expects the worst term to be the one where the fastly oscillating exponential is differentiated : 
\[ \int_1^t \int \int e^{i s \varphi} s \nabla_{\xi} \varphi m(\xi, \eta, \sigma) \widehat{f}(s, \xi-\eta-\sigma) \widehat{f}(s, \eta) \widehat{f}(s, \sigma) ~ d\sigma d\eta ds \]
In this context, we will search for a relation of the form : 
\begin{equation} m \nabla_{\xi} \varphi = \mu_1 \varphi + \mu_2 \nabla_{\eta} \varphi + \mu_3 \nabla_{\sigma} \varphi \label{equstar-heur} \end{equation}
which will allow to apply an integration by parts either in $s, \eta$ or $\sigma$ and hence go from a cubic to a quartic non-linearity in the first case (or a boundary in time term), or gain a $s^{-1}$ factor. 

The difficulty in the case of the mZK equation is that \eqref{equstar-heur} does not hold in general. However, one has that 
\[ \{ \varphi = 0, ~ \nabla_{\eta, \sigma} \varphi = 0 \} ~ \subset ~ \{ m \nabla_{\xi} \varphi = 0 \} \]
which means that $m \nabla_{\xi} \varphi$ vanishes on the resonant set. The idea is now to use \eqref{equstar-heur} whenever it holds, and near the singular points where \eqref{equstar-heur} does not hold, partition the space in order to use the cancellation of $m \nabla_{\xi} \varphi$ very close to the singular point, and use a singular version of \eqref{equstar-heur} away. This is very much in the spirit of the usual proof of the stationary phase lemma. 

\section{Kernel estimate} \label{section-dispersive-estimate}

Let us only consider the linear part of \eqref{equ_ZK_coordab} : 
\[ \partial_t u + \partial_{x_a}^3 u + \partial_{x_b}^3 u = 0 \]
The linear evolution is obtained by the kernel defined through its Fourier transform : 
\[ \widehat{K}_t(\xi_a, \xi_b) = e^{i t (\xi_a^3 + \xi_b^3)} \]
where $\xi_a$, $\xi_b$ are the Fourier dual variables associated to $x_a, x_b$. Then the solution of the linear problem with initial data $u_0$ is 
\[ u(t) = K_t \ast u_0 \]

Note that $\widehat{K}_t(\xi_a, \xi_b) = \widehat{k}_t(\xi_a) \widehat{k}_t(\xi_b)$ for $k_t$ the Airy fundamental solution : 
\[ \widehat{k}_t(\xi) = e^{i t \xi^3} \]
It is then well-known that 
\[ |k_t(x)| \lesssim t^{-\frac{1}{3}} \langle t^{-\frac{1}{3}} x \rangle^{-\frac{1}{4}} \]
(see eg \cite{LinaresPoncedispersivebook}, proposition 1.6) and therefore
\[ |K_t(x_a, x_b)| \lesssim t^{-\frac{2}{3}} \langle t^{-\frac{1}{3}} x_a \rangle^{-\frac{1}{4}} \langle t^{-\frac{1}{3}} x_b \rangle^{-\frac{1}{4}} \]

Likewise, for any $\beta \in \left[ 0, \frac{1}{2} \right]$, 
\[ \left| |\partial|^{\beta} k_t(x) \right| \lesssim t^{-\frac{1}{3}-\frac{\beta}{3}} \langle t^{-\frac{1}{3}} x \rangle^{-\frac{1}{4}+\frac{\beta}{2}} \]
and so for any $\beta_a, \beta_b \in \left[ 0, \frac{1}{2} \right]$, 
\[ \left| |\partial_a|^{\beta_a} |\partial_b|^{\beta_b} K_t(x_a, x_b) \right| \lesssim t^{-\frac{2}{3} - \frac{\beta_a+\beta_b}{3}} \langle t^{-\frac{1}{3}} x_a \rangle^{-\frac{1}{4}+\frac{\beta_a}{2}} \langle t^{-\frac{1}{3}} x_b \rangle^{-\frac{1}{4}+\frac{\beta_b}{2}} \]

From this we can deduce that $K_t$ lies in some weak Lebesgue spaces. Recall that the weak Lebesgue space $L^p_w(\mathbb{R}^d)$, defined for $p \in [1, \infty)$, is the set of functions $g$ such that 
\[ \Vert g \Vert_{L^p_w(\mathbb{R}^d)} := \sup_{\lambda > 0} \lambda^p \mathcal{L}^d\left( \left\{ z \in \mathbb{R}^d, ~ |g(z)| \geq \lambda \right\} \right) ~~ < \infty \]
where $\mathcal{L}^d$ is the $d$-dimensional Lebesgue measure, endowed with the topology induced by $\Vert \cdot \Vert_{L^p_w}$. Note that $\Vert \cdot \Vert_{L^p_w}$ is not a norm : we refer to \cite{SteinFourierbook}, chapter V, for a complete treatment of weak Lebesgue spaces and Lorentz spaces. 

\begin{Lem} Let $p \in [2, \infty)$ and $\beta_a, \beta_b \in \left[ 0, \frac{1}{2} \right)$ be such that 
\[ p \geq \max\left( \frac{4}{1-2\beta_a}, \frac{4}{1-2\beta_b} \right) \] 
Then $|\partial_a|^{\beta_a} |\partial_b|^{\beta_b} K_t \in L^p_w$ with 
\[ \Vert |\partial_a|^{\beta_a} |\partial_b|^{\beta_b} K_t \Vert_{L^p_w} \lesssim t^{-\frac{2+\beta_a+\beta_b}{3} + \frac{2}{3p}} \]
\label{lemKtweakLp}
\end{Lem}

\begin{Dem}
Let $p, \beta_a, \beta_b$ be as in the lemma. Then we have 
\[ \left| |\partial_a|^{\beta_a} |\partial_b|^{\beta_b} K_t(x_a, x_b) \right| \lesssim t^{-\frac{2}{3} - \frac{\beta_a+\beta_b}{3}} \langle t^{-\frac{1}{3}} x_a \rangle^{-\frac{1}{4}+\frac{\beta_a}{2}} \langle t^{-\frac{1}{3}} x_b \rangle^{-\frac{1}{4}+\frac{\beta_b}{2}} \]
But $\langle t^{-\frac{1}{3}} x_a \rangle^{-\frac{1}{4}+\frac{\beta_a}{2}}$ is in $L^q$ for any $q > p_a := \frac{4}{1 - 2 \beta_a}$ and in $L^{p_a}_w$. Likewise, $\langle t^{-\frac{1}{3}} x_b \rangle^{-\frac{1}{4}+\frac{\beta_b}{2}}$ is in $L^q$ for any $q > p_b := \frac{4}{1 - 2 \beta_b}$ and in $L^{p_b}_w$. Their norms in these spaces may be estimated by scaling as
\[ \Vert \langle t^{-\frac{1}{3}} x_a \rangle^{-\frac{1}{4}+\frac{\beta_a}{2}} \Vert_{L^q_w(\mathbb{R})} \leq C(q) t^{\frac{1}{3q}} \]
for any $q \in [p_a, \infty)$, where $C(q)$ is independent of $t$. In particular, we deduce that $K_t \in L^q_w$ for any $q \geq \max(p_a, p_b)$ with 
\[ \Vert K_t \Vert_{L^q_w(\mathbb{R}^2)} \leq C(q) t^{\frac{2}{3q}-\frac{2}{3} - \frac{\beta_a+\beta_b}{3}} \]
and the hypothesis on $p$ allows to choose $q = p$. 
\end{Dem}

Recall the following classical theorem : 

\begin{Theo}[Hardy-Littlewood-Sobolev] Let $d \in \mathbb{N}$, $p, q, r \in [1, \infty]$ be such that $\frac{1}{p} + \frac{1}{q} = \frac{1}{r} + 1$. There exists $C = C(p, q, r, d) > 0$ such that for any $F \in L^p(\mathbb{R}^d)$, $G \in L^q_w(\mathbb{R}^d)$, then $F \ast G \in L^r(\mathbb{R}^d)$ and 
\[ \Vert F \ast G \Vert_{L^r} \leq C \Vert F \Vert_{L^p} \Vert G \Vert_{L^q_w} \]
\label{theoHLS}
\end{Theo}

See for instance \cite{SteinFourierbook}. 

\section{Multilinear estimates and paraproducts} \label{section-paraproducts}

Recall the following classical linear estimate : 

\begin{Theo}[Hörmander-Mikhlin] Let $m : \mathbb{R}^2 \to \mathbb{C}$ be a bounded function, such that 
\[ \Vert m \Vert_{HM} := \sup_{\substack{\alpha \in \mathbb{N}^2, \\|\alpha| \leq d+1}} \sup_{\xi \in \mathbb{R}^2} |\xi|^{|\alpha|} |\partial^{\alpha} m(\xi)| ~~ < \infty \]
Then the linear operator 
\[ m(D) : ~ f \in \mathcal{S}(\mathbb{R}^2) ~ \mapsto ~ \mathcal{F}^{-1} \left[ m(\xi) \widehat{f}(\xi) \right] \in \mathcal{S}'(\mathbb{R}^2) \]
extends to a continuous linear operator $L^p(\mathbb{R}^2) \to L^p(\mathbb{R}^2)$ for any $1 < p < \infty$, with 
\[ \Vert m(D) \Vert_{L^p \to L^p} \leq C(p) \Vert m \Vert_{HM} \]
for some universal constant $C(p)$ depending only on $p$. \label{theoHM}
\end{Theo}

In the multilinear case, we have the following extension : 

\begin{Theo}[Coifman-Meyer \cite{CoifmanMeyer}] Let $n \geq 2$ and $m : \mathbb{R}^{2n} \to \mathbb{C}$ be a bounded function, such that 
\[ \Vert m \Vert_{CM} := \sup_{\alpha \in \mathbb{N}^{2n}, |\alpha| \leq N} \sup_{(\eta_1, ..., \eta_n) \in \mathbb{R}^{2n}} |(\eta_1, ..., \eta_n)|^{|\alpha|} |\partial^{\alpha} m(\eta_1, ..., \eta_n)| ~~ < \infty \]
where $N \in \mathbb{N}$ is a large enough (universal) number. Then the $n$-linear operator 
\[ T_m : ~ (f_1, ..., f_n) \in \mathcal{S}(\mathbb{R}^2)^n ~ \mapsto ~ \mathcal{F}^{-1}_{\xi \to x} \left[ \int m(\eta) \prod_{i = 1}^n \widehat{f_i}(\eta_i) ~ \delta(\xi-\eta_1-...-\eta_n) d\eta_1 ... d\eta_n \right] \in \mathcal{S}'(\mathbb{R}^2) \]
extends to a continuous $n$-linear operator $L^{p_1}(\mathbb{R}^2) \times ... \times L^{p_n}(\mathbb{R}^2) \to L^q(\mathbb{R}^2)$ for any $1 < p_i \leq \infty$, $0 < q < \infty$ such that $\frac{1}{q} = \sum_{i = 1}^n \frac{1}{p_i}$, with 
\[ \Vert T_m \Vert_{L^{p_1} \times ... \times L^{p_n} \to L^q} \leq C(p_1, ..., p_n, q) \Vert m \Vert_{CM} \]
for some universal constant $C(p, q, r)$ depending only on $q$ and the $p_i$. 
\label{theoCM}
\end{Theo}

Note that, for $m \equiv 1$, one recovers the Hölder inequality (without the optimal constant). For a proof of these results, we refer the reader to the books \cite{MuscaluHAbook1, MuscaluHAbook2}. 

We will need the following extension that allows to treat separately the two coordinates :  

\begin{Theo}[Muscalu-Pipher-Tao-Thiele, \cite{MPTTBiparameter2004}] Let $m : \mathbb{R}^4 \to \mathbb{C}$ be a bounded function, such that 
\[ \Vert m \Vert_{MPTT} := \sup_{\substack{\alpha \in \mathbb{N}^2, \\ \beta \in \mathbb{N}^2, \\ |\alpha| + |\beta| \leq N}} ~ \sup_{\substack{(\xi, \eta) = (\xi_a, \xi_b, \eta_a, \eta_b) \\ \in \mathbb{R}^2 \times \mathbb{R}^2}} |(\xi_a, \eta_a)|^{|\alpha|} |(\xi_b, \eta_b)|^{|\beta|} |\partial_{\xi_a}^{\alpha_1} \partial_{\xi_b}^{\beta_1} \partial_{\eta_a}^{\alpha_2} \partial_{\eta_b}^{\beta_2} m(\xi, \eta)| ~~ < \infty \]
where $N \in \mathbb{N}$ is a large enough (universal) number. Then the bilinear operator 
\[ T_m : ~ (f, g) \in \mathcal{S}(\mathbb{R}^2)^2 ~ \mapsto ~ \mathcal{F}^{-1} \left[ \int m(\xi, \eta) \widehat{f}(\xi-\eta) \widehat{g}(\eta) ~ d\eta \right] \in \mathcal{S}'(\mathbb{R}^2) \]
extends to a continuous bilinear operator $L^p(\mathbb{R}^2) \times L^q(\mathbb{R}^2) \to L^r(\mathbb{R}^2)$ for any $1 < p, q \leq \infty$, $0 < r < \infty$ such that $\frac{1}{r} = \frac{1}{p} + \frac{1}{q}$, with 
\[ \Vert T_m \Vert_{L^p \times L^q \to L^r} \leq C(p, q, r) \Vert m \Vert_{MPTT} \]
for some universal constant $C(p, q, r)$ depending only on $(p, q, r)$. 

Similarly, let $m' : \mathbb{R}^6 \to \mathbb{C}$ be a bounded function, such that 
\begin{align*} 
\Vert m \Vert_{MPTT} := \sup \Big\{ &|(\xi_a, \eta_a, \sigma_a)|^{|\alpha|} |(\xi_b, \eta_b, \sigma_b)|^{|\beta|} |\partial_{\xi_a}^{\alpha_1} \partial_{\xi_b}^{\beta_1} \partial_{\eta_a}^{\alpha_2} \partial_{\eta_b}^{\beta_2} \partial_{\sigma_a}^{\alpha_3} \partial_{\sigma_b}^{\beta_3} m(\xi, \eta, \sigma)|, \\
&\quad \alpha \in \mathbb{N}^3, \beta \in \mathbb{N}^3, |\alpha| + |\beta| \leq N, (\xi, \eta, \sigma) \in \mathbb{R}^6, \Big\} ~~ < \infty 
\end{align*}
where $N$ is a large enough universal number. Then the trilinear operator 
\[ T_m : ~ (f, g, h) \in \mathcal{S}(\mathbb{R}^2)^3 ~ \mapsto ~ \mathcal{F}^{-1} \left[ \int m(\xi, \eta, \sigma) \widehat{f}(\xi-\eta-\sigma) \widehat{g}(\eta) \widehat{h}(\sigma) ~ d\eta d\sigma \right] \in \mathcal{S}'(\mathbb{R}^2) \]
extends to a continuous trilinear operator $L^p(\mathbb{R}^2) \times L^q(\mathbb{R}^2) \times L^r(\mathbb{R}^2) \to L^s(\mathbb{R}^2)$ for any $1 < p, q, r \leq \infty$, $0 < s < \infty$ such that $\frac{1}{s} = \frac{1}{p} + \frac{1}{q} + \frac{1}{r}$, with 
\[ \Vert T_m \Vert_{L^p \times L^q \times L^r \to L^q} \leq C(p, q, r, s) \Vert m \Vert_{MPTT} \]
for some universal constant $C(p, q, r, s)$ depending only on $(p, q, r, s)$. 
\label{theobiparaproduct}
\end{Theo}

\begin{Rem} We restrict our attention to the bilinear and trilinear cases, but \cite{MPTTBiparameter2004} covers any $n$-linear setting. 
\end{Rem}

\section{A priori estimate}

\subsection{Statement of the a priori estimate} \label{section-statement-apriori}

Let $u : [1, T] \times \mathbb{R}^2 \to \mathbb{R}$ be a solution of the Zakharov-Kuznetsov equation \eqref{equ_ZK_coordab}, for some initial data $u(1)$ givent at time $t = 1$, and $T > 1$ be an existence time. Define the profile function as 
\[ f(t) = K_{-t} \ast u(t) \]
where $K_t = e^{tL}$ is defined as the kernel of the linear evolution (see section \ref{section-dispersive-estimate}) and $L$ the linear operator of \eqref{equ_ZK_coordab}, so that $f$ satisfies
\[ \partial_t f = e^{tL} \left( \partial_{x_a} + \partial_{x_b} \right) \left( e^{-tL} f \right)^3 \]

Let us define the norm
\[ \Vert u \Vert_{X, T} = \sup_{t \in [1, T]} \left( \Vert u(t) \Vert_{H^3}, ~ \Vert x_a f(t) \Vert_{L^2}, ~ \Vert x_b f(t) \Vert_{L^2} \right) \]

We aim to prove : 

\begin{Prop} If $u$ satisfies \eqref{equ_ZK_coordab}, then 
\[ \Vert u \Vert_{X, T} \leq \Vert u \Vert_{X, 1} + C \Vert u \Vert_{X, T}^2 \left( 1 + \Vert u \Vert_{X, T}^4 \right) \]
for some universal constant $C > 0$. \label{prop-est-a-priori}
\end{Prop}

\subsection{Preliminary estimates} \label{section-prelest}

\begin{Lem}[dispersive estimate] Let $\varepsilon > 0$. Then for any $t \in [1, T]$, 
\begin{align*}
\Vert u(t) \Vert_{L^{\infty}} &\lesssim t^{-\frac{2}{3}+\varepsilon} \Vert u \Vert_{X, T} \\
\Vert \nabla u(t) \Vert_{L^{\infty}} &\lesssim t^{-\frac{5}{9}+\varepsilon} \Vert u \Vert_{X, T} \\
\Vert |\partial_{x_a}|^{\frac{1}{2}} u(t) \Vert_{L^{\infty}} &\lesssim t^{-\frac{2}{3}-\frac{1}{6}+\varepsilon} \Vert u \Vert_{X, T} \\
\Vert |\partial_{x_b}|^{\frac{1}{2}} u(t) \Vert_{L^{\infty}} &\lesssim t^{-\frac{2}{3}-\frac{1}{6}+\varepsilon} \Vert u \Vert_{X, T} 
\end{align*} \label{lemLinfest}
\end{Lem}

\begin{Dem}
Note that, without loss of generality, we may assume that $\varepsilon > 0$ is small. 

Recall that $u(t) = K_t \ast f(t)$, and by lemma \ref{lemKtweakLp}, $K_t \in L^p_w(\mathbb{R}^2)$ for any $p \geq 4$. In particular, choosing $p = \frac{2}{3 \varepsilon}$, which is larger than $4$ if $\varepsilon$ is small enough, we get by lemma \ref{lemKtweakLp}, theorem \ref{theoHLS} and Hölder's inequality that 
\begin{align*}
\Vert u(t) \Vert_{L^{\infty}} &\lesssim \Vert f(t) \Vert_{L^{p'}} \Vert K_t \Vert_{L^p_w} \\
&\lesssim \Vert (1 + |x_a| + |x_b|) f(t) \Vert_{L^2} t^{-\frac{2}{3} + \frac{2}{3p}} \\
&\lesssim t^{-\frac{2}{3}+\varepsilon} \left( \Vert u(t) \Vert_{L^2} + \Vert (|x_a| + |x_b|) f(t) \Vert_{L^2} \right) \\
&\lesssim t^{-\frac{2}{3}+\varepsilon} \Vert u \Vert_{X, T}
\end{align*}
The third and fourth estimates are done the same way : by lemma \ref{lemKtweakLp}, $|\partial_{x_a}|^{\frac{1}{2}-\varepsilon} K_t \in L^p_w(\mathbb{R}^2)$ for $p = \frac{2}{\varepsilon}$, and so by theorem \ref{theoHLS}, 
\begin{align*}
\Vert |\partial_{x_a}|^{\frac{1}{2}-\varepsilon} u(t) \Vert_{L^{\infty}} &\lesssim \Vert f(t) \Vert_{L^{p'}} \Vert |\partial_{x_a}|^{\frac{1}{2}-\varepsilon} K_t \Vert_{L^p_w} \\
&\lesssim \Vert u \Vert_{X, T} t^{-\frac{2}{3}-\frac{1}{6}+\frac{2 \varepsilon}{3} + \frac{\varepsilon}{3}} = t^{-\frac{2}{3}-\frac{1}{6}+\varepsilon} \Vert u \Vert_{X, T} 
\end{align*}
To get an estimate on $\Vert |\partial_{x_a}|^{\frac{1}{2}} u(t) \Vert_{L^{\infty}}$, we write by Sobolev's embedding : 
\[ \Vert |\partial_{x_a}|^{\frac{1}{2}} u(t) \Vert_{L^{\infty}} \lesssim \Vert |\partial_{x_a}|^{\frac{1}{2}} u(t) \Vert_{W^{3\delta, \frac{1}{\delta}}} \]
for any parameter $\delta > 0$ that we will choose small. Then by interpolation : 
\[ \Vert |\partial_{x_a}|^{\frac{1}{2}} u(t) \Vert_{L^{\frac{1}{\delta}}} \lesssim \Vert |\partial_{x_a}|^{\frac{1}{2}-\delta} u(t) \Vert_{L^{\frac{1}{\delta}}}^{\frac{1}{1+\delta}} \Vert \partial_{x_a} u(t) \Vert_{L^{\frac{1}{\delta}}}^{\frac{\delta}{1+\delta}} \]
Now $\Vert \partial_{x_a} u(t) \Vert_{L^{\frac{1}{\delta}}} \lesssim \Vert u(t) \Vert_{H^3} \leq \Vert u \Vert_{X, T}$ by Sobolev's embedding, while by interpolation again, 
\[ \Vert |\partial_{x_a}|^{\frac{1}{2}-\delta} u(t) \Vert_{L^{\frac{1}{\delta}}} \leq \Vert |\partial_{x_a}|^{\frac{1}{2}-\delta} u(t) \Vert_{L^2}^{2\delta} \Vert |\partial_{x_a}|^{\frac{1}{2}-\delta} u(t) \Vert_{L^{\infty}}^{1-2\delta} \]
A similar estimate can be proved by interpolation on $\Vert |\partial_{x_a}|^{\frac{1}{2}} u(t) \Vert_{\dot{W}^{3\delta, \frac{1}{\delta}}}$, so we conclude that 
\[ \Vert |\partial_{x_a}|^{\frac{1}{2}} u(t) \Vert_{L^{\infty}} \lesssim \Vert |\partial_{x_a}|^{\frac{1}{2}-\delta} u(t) \Vert_{L^{\infty}}^{\theta} \Vert u \Vert_{X, T}^{1-\theta} \]
for some $\theta = \theta(\delta)$ that converges to $1$ as $\delta \to 0$. In particular, by choosing $\delta$ small enough with respect to $\varepsilon$, we obtain 
\[ \Vert |\partial_{x_a}|^{\frac{1}{2}} u(t) \Vert_{L^{\infty}} \lesssim t^{-\frac{2}{3}-\frac{1}{6}+\varepsilon} \Vert u \Vert_{X, T} \]
The estimate with $\partial_{x_b}$ is symmetric. 

For the second estimate, we prove it replacing $\nabla$ by $\partial_{x_a}$ and conclude by symmetry of $x_a$ and $x_b$. Then, for a small parameter $\delta > 0$, we have by Sobolev's embedding : 
\[ \Vert \partial_{x_a} u(t) \Vert_{L^{\infty}} \lesssim \Vert \partial_{x_a} u(t) \Vert_{W^{3 \delta, \frac{1}{\delta}}} \]
Now we may interpolate : 
\[ \Vert \partial_{x_a} u(t) \Vert_{L^{\frac{1}{\delta}}} \lesssim \Vert |\partial_{x_a}|^{\frac{1}{2}} u(t) \Vert_{L^{\frac{1}{\delta}}}^{\frac{2}{3}} \Vert \partial_{x_a}^2 u(t) \Vert_{L^{\frac{1}{\delta}}}^{\frac{1}{3}} \]
But by Sobolev's embedding, 
\[ \Vert \partial_{x_a}^2 u(t) \Vert_{L^{\frac{1}{\delta}}} \lesssim \Vert u(t) \Vert_{H^3} \lesssim \Vert u \Vert_{X, T} \]
On the other hand, by interpolation again, 
\[ \Vert |\partial_{x_a}|^{\frac{1}{2}} u(t) \Vert_{L^{\frac{1}{\delta}}} \leq \Vert |\partial_{x_a}|^{\frac{1}{2}} u(t) \Vert_{L^2}^{2\delta} \Vert |\partial_{x_a}|^{\frac{1}{2}} u(t) \Vert_{L^{\infty}}^{1-2\delta} \]
and we already proved an estimate on these factors. This leads to 
\[ \Vert \partial_{x_a} u(t) \Vert_{L^{\frac{1}{\delta}}} \lesssim \Vert u \Vert_{X, T} \left( t^{-\frac{2}{3}-\frac{1}{6}+\delta} \right)^{\frac{2}{3} (1 - 2 \delta)} = t^{-\kappa} \Vert u \Vert_{X, T} \]
where $\kappa = \kappa(\delta)$ and converges to $\frac{2}{3} \left( \frac{2}{3} + \frac{1}{6} \right) = \frac{5}{9}$ when $\delta \to 0$. A similar bound can be proved the same way on $\Vert \partial_{x_a} u(t) \Vert_{\dot{W}^{3\delta, \frac{1}{\delta}}}$, and so by choosing $\delta$ small enough with respect to $\varepsilon$, we get the desired result. 
\end{Dem}

\subsection{Proof of theorem \ref{main-theo}} \label{section-proof-mainth}

At this point, we can prove already theorem \ref{main-theo}, assuming proposition \ref{prop-est-a-priori} holds. 

We start by the following lemma : 

\begin{Lem} For any $1 \leq t \leq T$, 
\[ \Vert \partial_t f(t) \Vert_{H^2} \lesssim t^{-1-\frac{1}{9}} \Vert u \Vert_X^3 \] \label{lemcontrdtf}
\end{Lem}

\begin{Dem}
We have that
\[ \partial_t f(t) = e^{tL} \left( \partial_{x_a} + \partial_{x_b} \right) \left( u(t) \right)^3 \]
so that for any $\delta > 0$ small, 
\begin{align*}
\Vert \partial_t f(t) \Vert_{H^2} &\lesssim \Vert u(t)^3 \Vert_{H^3} \\
&\lesssim \Vert u(t) \Vert_{L^{\infty}}^2 \Vert u(t) \Vert_{H^3} + \Vert u(t) \Vert_{L^{\infty}} \Vert u(t) \Vert_{W^{1, \infty}} \Vert u(t) \Vert_{H^2} + \Vert u(t) \Vert_{W^{1, 6}}^3 \\
&\lesssim \Vert u \Vert_{X, T}^3 t^{-\frac{4}{3}+2\delta} + \Vert u \Vert_{X, T}^3 t^{-\frac{2}{3}-\frac{5}{9}+2\delta} + \left( \Vert u(t) \Vert_{H^3}^{1/3} \Vert u(t) \Vert_{L^{\infty}}^{2/3} \right)^3 \\
&\lesssim \Vert u \Vert_{X, T}^3 t^{-\frac{11}{9}+2\delta}
\end{align*}
Choosing $\delta > 0$ small enough, we conclude. 
\end{Dem}

Now we prove theorem \ref{main-theo}. 

\begin{Dem}[of theorem \ref{main-theo}]
Note first that the hypothesis of theorem \ref{main-theo} control $\Vert u \Vert_{X, 1}$ : $\Vert u(1) \Vert_{H^3}$ is directly controlled by hypothesis and 
\begin{align*} 
\Vert x_a f(1) \Vert_{L^2} &= \Vert x_a e^{-L} u(1) \Vert_{L^2} \\
&\leq \Vert 3 \partial_{x_a}^2 e^{-L} u(1) \Vert_{L^2} + \Vert e^{-L} (x_a u(1)) \Vert_{L^2} \\
&\lesssim \Vert u(1) \Vert_{H^2} + \Vert x_a u(1) \Vert_{L^2}
\end{align*}
Likewise for $\Vert x_b f(1) \Vert_{L^2}$. 

Let $\varepsilon = \Vert u \Vert_{X, 1}$ be the size of the initial data, with $\varepsilon \leq 1$. 

Let now $T$ be an existence time such that $\Vert u \Vert_{X, T} \leq 2 \varepsilon$. Then one has by proposition \ref{prop-est-a-priori} : 
\begin{align*} 
\Vert u \Vert_{X, T} &\leq \Vert u \Vert_{X, 1} + C \Vert u \Vert_{X, T}^2 \left( 1 + \Vert u \Vert_{X, T}^4 \right) \\
&\leq \varepsilon + 4 C \varepsilon^2 \left( 1 + 16 \varepsilon^4 \right) \\
&\leq \varepsilon \left( 1 + 4C \varepsilon * 17 \right) 
\end{align*}
Now if $\varepsilon < \frac{1}{8*17C}$, then 
\[ \Vert u \Vert_{X, T} \leq \frac{3}{2} \varepsilon < 2 \varepsilon \]
By continuity of $T \mapsto \Vert u \Vert_{X, T}$, we deduce that we may choose $T$ arbitrarily large, and that the solution $u(t)$ remains in $H^3$ for all times. 

Now we just proved that $\Vert u \Vert_{X, T} \leq 2 \varepsilon$ for any $T$ as soon as $\varepsilon$ is small enough, and by lemma \ref{lemcontrdtf}, we get that $\partial_t f(t) \in L^1_t H^2_{x_a, x_b}$, which ensures that $f(t)$ converges in $H^2$. (Moreover, the limit $f(\infty)$ is $\varepsilon^3$-close of the initial data $f(1)$.) 
\end{Dem} 

In the rest of this paper, we simply write $\Vert \cdot \Vert_X$ and forget about the dependence in $T$. 

\section{Energy estimate} \label{section-H3est}

In this section, we prove the following estimate : 

\begin{Prop} There exists a universal constant $C > 0$ such that, for any $1 \leq t \leq T$, 
\[ \Vert u(t) \Vert_{H^3} \leq \Vert u(1) \Vert_{H^3} + C \Vert u \Vert_X^2 \]
\end{Prop}

The $L^2$-norm of $u$ is preserved by the flow, so we immediately get
\[ \Vert u(t) \Vert_{L^2} = \Vert u(1) \Vert_{L^2} \leq \Vert u(1) \Vert_{L^2} + C \Vert u \Vert_X^2 \]

On the other hand, 
\begin{align*}
\partial_t \Vert \partial_{x_a}^3 f(t) \Vert_{L^2}^2 
&= 2 \int \partial_{x_a}^3 f(t, x_a, x_b) \partial_{x_a}^3 \partial_t f(t, x_a, x_b) ~ dx_a dx_b \\
&= 2 \int \partial_{x_a}^3 u(t, x_a, x_b) \partial_{x_a}^3 \left( \partial_{x_a} + \partial_{x_b} \right) \left( u^3(t, x_a, x_b) \right) ~ dx_a dx_b \\
&= 6 \int \partial_{x_a}^3 u(t, x_a, x_b) \left( (\partial_{x_a} + \partial_{x_b}) \partial_{x_a}^3 u(t, x_a, x_b) \right) ~ u^2(t, x_a, x_b) ~ dx_a dx_b \\
&\quad + 6 \int \partial_{x_a}^3 u(t, x_a, x_b) ~ \left[ \partial_{x_a}^3, u^2(t, x_a, x_b) (\partial_{x_a} + \partial_{x_b}) \right] u(t, x_a, x_b) ~ dx_a dx_b \\
&= I + II
\end{align*}

The first term can be integrated by parts : 
\begin{align*}
I &= 6 \int \partial_{x_a}^3 u(t, x_a, x_b) \left( (\partial_{x_a} + \partial_{x_b}) \partial_{x_a}^3 u(t, x_a, x_b) \right) ~ u(t, x_a, x_b)^2 ~ dx_a dx_b \\
&= 6 \int (\partial_{x_a} + \partial_{x_b}) \left( (\partial_{x_a}^3 u(t, x_a, x_b))^2 \right) ~ u(t, x_a, x_b)^2 ~ dx_a dx_b \\
&= - 12 \int (\partial_{x_a}^3 u(t, x_a, x_b))^2 ~ u(t, x_a, x_b) (\partial_{x_a} + \partial_{x_b}) u(t, x_a, x_b) ~ dx_a dx_b 
\end{align*}
and then estimated by Hölder's inequality : 
\[ |I| \lesssim \Vert u(t) \Vert_{H^3}^2 \Vert \nabla u(t) \Vert_{L^{\infty}} \Vert u(t) \Vert_{L^{\infty}} \]

The second term can be estimated by Hölder's inequality and Kato-Ponce's commutator estimate \cite{KatoPonceCommut} : 
\begin{align*}
|II| &\lesssim \Vert u(t) \Vert_{H^3} \left\Vert \left[ \partial_{x_a}^3, u^2(t) (\partial_{x_a} + \partial_{x_b}) \right] u(t) \right\Vert_{L^2} \\
&\lesssim \Vert u(t) \Vert_{H^3} \left( \Vert \nabla (u^2(t)) \Vert_{L^{\infty}} \Vert u(t) \Vert_{H^3} + \Vert u^2(t) \Vert_{H^3} \Vert \nabla u(t) \Vert_{L^{\infty}} \right) \\
&\lesssim \Vert u(t) \Vert_{H^3}^2 \Vert u(t) \Vert_{L^{\infty}} \Vert u(t) \Vert_{L^{\infty}} 
\end{align*}

Putting everything together, we deduce by lemma \ref{lemLinfest} : 
\begin{align*}
\partial_t \Vert \partial_{x_a}^3 f(t) \Vert_{L^2}^2 
&\lesssim \Vert u(t) \Vert_{H^3}^2 \Vert \nabla u(t) \Vert_{L^{\infty}} \Vert u(t) \Vert_{L^{\infty}} \\
&\lesssim \Vert u \Vert_X^4 t^{-\frac{11}{9}+2\delta}
\end{align*}
for any $\delta > 0$.Choosing $\delta$ small enough, we deduce that for any $1 \leq t \leq T$, 
\[ \Vert \partial_{x_a}^3 f(t) \Vert_{L^2}^2 \leq \Vert \partial_{x_a}^3 f(1) \Vert_{L^2}^2 + C \Vert u \Vert_X^4 \]
for some constant $C > 0$ universal. This implies
\[ \Vert \partial_{x_a}^3 u(t) \Vert_{L^2} \leq \Vert \partial_{x_a}^3 u(1) \Vert_{L^2} + C \Vert u \Vert_X^2 \]
(for another universal constant $C$). 

By symmetry between $x_a$ and $x_b$, we conclude that 
\[ \Vert u(t) \Vert_{H^3} \leq \Vert u(1) \Vert_{H^3} + C \Vert u \Vert_X^2 \]

\section{Weighted estimate} \label{section-weightedest}

In this section, we prove the following estimate : 

\begin{Prop} There exists $C > 0$ universal such that, for any $1 \leq t \leq T$, 
\[ \Vert x_a f(t) \Vert_{L^2} \leq \Vert x_a f(1) \Vert_{L^2} + C \Vert u \Vert_X^3 \left( 1 + \Vert u \Vert_X^2 \right) \] \label{propaprioripoids}
\end{Prop}

By symmetry between $x_a$ and $x_b$, proposition \ref{prop-est-a-priori} will follow. 

To prove such an estimate, we write using Duhamel formula and applying a Fourier transform : 
\begin{align*}
\widehat{f}(t, \xi) &= \widehat{f}(1, \xi) + \int_1^t \int \int e^{i s \varphi(\xi, \eta, \sigma)} (\xi_a + \xi_b) \widehat{f}(s, \eta) \widehat{f}(s, \sigma) \widehat{f}(s, \xi-\eta-\sigma) ~ d\eta d\sigma ds 
\end{align*}
where $\xi = (\xi_a, \xi_b) \in \mathbb{R}^2$, and likewise for $\eta, \sigma$, and 
\[ \varphi(\xi, \eta, \sigma) = \xi_a^3 + \xi_b^3 - \eta_a^3 - \eta_b^3 - \sigma_a^3 - \sigma_b^3 - (\xi_a-\eta_a-\sigma_a)^3 - (\xi_b-\eta_b-\sigma_b)^3 \]
is the interaction phase. 

We will denote by $\rho = \xi - \eta - \sigma$. 

Then the weighted estimate will follow from an estimate in $L^2$ of 
\begin{equation} \partial_{\xi_a} \int_1^t \int \int e^{i s \varphi(\xi, \eta, \sigma)} (\xi_a + \xi_b) \widehat{f}(s, \eta) \widehat{f}(s, \sigma) \widehat{f}(s, \rho) ~ d\eta d\sigma \label{termtotpoids} \end{equation}

\subsection{Space-time resonances analysis} \label{section-spacetimeresonance}

The space-time resonances is based on the following idea. Considering a term in the form of \eqref{termtotpoids}, one applies a stationary phase heuristics to deduce that the main contributions are coming from points where the phase $s \varphi$ is stationary. Given that there is an integration in both frequency $\eta, \sigma$ and time $s$, we infer that the behaviour of the solution is driven by its behaviour near points $(\xi, \eta, \sigma)$ such that 
\[ \nabla_{\eta, \sigma} \varphi(\xi, \eta, \sigma) = 0, \quad \varphi(\xi, \eta, \sigma) = 0 \]

\begin{Def} The set $\mathcal{S} = \{ (\xi, \eta, \sigma) \in \mathbb{R}^6, ~ \nabla_{\eta, \sigma} \varphi(\xi, \eta, \sigma) = 0 \}$ is called the set of space resonances. 

The set $\mathcal{T} = \{ (\xi, \eta, \sigma) \in \mathbb{R}^6, ~ \varphi(\xi, \eta, \sigma) = 0 \}$ is called the set of time resonances. 

The set $\mathcal{R} = \mathcal{S} \cap \mathcal{T}$ is called the set of space-time resonances. 
\end{Def}

\begin{Lem} The space resonant points satisfy 
\[ |\eta_a| = |\sigma_a| = |\rho_a|, \quad |\eta_b| = |\sigma_b| = |\rho_b| \]
Among space resonant points, if we denote $\eta_a = \epsilon_1^a \rho_a, \sigma_a = \epsilon_2^a \rho_a$ and $\eta_b = \epsilon_1^b \rho_b, \sigma_b = \epsilon_2^b \rho_b$ where $\epsilon \in \{ -1, 1 \}$, then the space-time resonant points also satisfy $\Big[ \big( \epsilon_1^a = -1$ or $\epsilon_2^a = -1$ or $\rho_a = 0 \big)$ and $\big( \epsilon_1^b = -1$ or $\epsilon_2^b = -1$ or $\rho_b = 0 \big) \Big]$ or $\Big[ \epsilon_1^a = \epsilon_2^a = \epsilon_1^b = \epsilon_2^b = 1$ and $\rho_a + \rho_b = 0 \Big]$. \label{carac_res_esptps}
\end{Lem}

\begin{Dem}
We compute : 
\[ \partial_{\eta_a} \varphi(\xi, \eta, \sigma) = 3 (\rho_a^2 - \eta_a^2), \quad \partial_{\sigma_a} \varphi(\xi, \eta, \sigma) = 3 (\rho_a^2 - \sigma_a^2) \]
This shows the first statement by symmetry between the $a$ and $b$ coordinates. 

For space-time resonances, we write given such a space resonant point : 
\[ \begin{aligned}
\varphi(\xi, \eta, \sigma) &= (\eta_a + \rho_a + \sigma_a)^3 - \eta_a^3 - \sigma_a^3 - \rho_a^3 + (\eta_b + \rho_b + \sigma_b)^3 - \eta_b^3 - \sigma_b^3 - \rho_b^3 \\
&= \rho_a^3 \left( (1 + \epsilon_1^a + \epsilon_2^a)^3 - \epsilon_1^a - \epsilon_2^a - 1 \right) + \rho_b^3 \left( (1 + \epsilon_1^b + \epsilon_2^b)^3 - \epsilon_1^b - \epsilon_2^b - 1 \right) 
\end{aligned} \]
But $1 + \epsilon_1^a + \epsilon_2^a \in \{ -1, 1, 3 \}$. In the first two cases, the term in front of $\rho_a^3$ vanishes, which is equivalent to $\epsilon_1^a + \epsilon_2^a \leq 1$ ; in the third case, it is equal to $24$. We may therefore check that the only possibility for $\varphi$ to vanish as well are exhausted by the second statement of the lemma. 
\end{Dem}

By the above computation of the space-time resonant set, we may now try to obtain structure on the non-linear term. The most favorable case happens when the symbol of the non-linearity, here $\xi_a + \xi_b$, can be expressed as a combination of $\varphi$ and $\nabla_{\eta, \sigma} \varphi$~: that way, one can apply integrations by parts in $s$, $\eta$ or $\sigma$ and recover terms that are easier to control. 

In our case, such a relation does not hold, but considering the presence of a $\partial_{\xi_a}$, we will try to obtain a relation of the form~: 
\[ (\xi_a + \xi_b) \partial_{\xi_a} \varphi(\xi, \eta, \sigma) = \mu_0(\xi, \eta, \sigma) \varphi(\xi, \eta, \sigma) + \mu_1(\xi, \eta, \sigma) \cdot \nabla_{\eta, \sigma} \varphi(\xi, \eta, \sigma) \]
for $\mu_0, \mu_1$ symbols in some class. 

\begin{Lem} There exists symbols $m_1^a, m_1^b, m_2^a, m_2^b, m_3$ with a finite $\Vert \cdot \Vert_{MPTT}$-norm, such that, if we denote by $\omega = (\xi, \eta, \sigma, \rho)$ and $\omega_a = (\xi_a, \eta, \sigma_a, \rho_a)$, then 
\begin{equation} \begin{aligned}
(\xi_a + \xi_b) \partial_{\xi_a} \varphi(\xi, \eta, \sigma) &= \partial_{\eta_a} \varphi(\xi, \eta, \sigma) m_1^a(\xi, \eta, \sigma) \cdot \omega + \partial_{\sigma_a} \varphi(\xi, \eta, \sigma) m_2^a(\xi, \eta, \sigma) \cdot \omega \\
&+ \frac{\partial_{\eta_b} \varphi(\xi, \eta, \sigma)}{|\omega_a|} m_1^b(\xi, \eta, \sigma) : \omega^2 
+ \frac{\partial_{\sigma_b} \varphi(\xi, \eta, \sigma)}{|\omega_a|} m_2^b(\xi, \eta, \sigma) : \omega^2 \\
&+ \frac{\varphi(\xi, \eta, \sigma) }{|\omega_a|} m_3(\xi, \eta, \sigma) \cdot \omega
\end{aligned} \label{rel_const} \end{equation}
Moreover, the singularities in $\frac{1}{|\omega_a|}$ only appear near the space-time resonant points. \label{lem_rel_const}
\end{Lem}

\begin{Dem}
We first notice that, if we prove the above relation on the sphere $\left\{ \eta_a^2 + \eta_b^2 + \sigma_a^2 + \sigma_b^2 + \rho_a^2 + \rho_b^2 = 1 \right\}$, then the same relation can be extended to the whole space by homogeneity. 

We first prove the result near any point of the sphere. 

Let $Z = (\eta^Z, \sigma^Z, \rho^Z)$ be on this sphere. If $Z$ is a space-time resonant point, we set $\epsilon_i^{\alpha}$, $i = 1, 2$, $\alpha = a, b$ as in lemma \ref{carac_res_esptps}. 
\begin{itemize}[label=\textbullet]
\item If $Z$ is not a space-time resonant point, then one of the functions $\nabla_{\eta} \varphi$, $\nabla_{\sigma} \varphi$, $\varphi$ does not vanish in a neighborhood of $Z$, and we may construct locally these symbols without any singularity. 
\item If $Z$ is a space-time resonant point such that $\rho_a^Z \neq 0$ and $\epsilon_1^a + \epsilon_2^a + 1 \in \{ -1, 1 \}$, we have that, near $Z$, 
\[ \partial_{\xi_a} \varphi(\xi, \eta) = 3 (\xi_a^2 - \rho_a^2) \]
vanishes on $\{ |\eta_a| = |\sigma_a| = |\rho_a| \} = \{ \partial_{\eta_a} \varphi = \partial_{\sigma_a} \varphi = 0 \}$. However, 
\[ \nabla \partial_{\eta_a} \varphi = 3 \begin{pmatrix} - 2 \eta_a \\ 0 \\ 0 \\ 0 \\ 2 \rho_a \\ 0 \end{pmatrix}, \quad \nabla \partial_{\sigma_a} \varphi = 3 \begin{pmatrix} 0 \\ 0 \\ - 2 \sigma_a \\ 0 \\ 2 \rho_a \\ 0 \end{pmatrix} \]
(where $\nabla$ differentiates with respect to coordinates $(\eta, \sigma, \rho)$) and these two vectors are linearly independent, so $\partial_{\eta_a} \varphi$ and $\partial_{\sigma_a} \varphi$ are local normal coordinates to the local manifold $\{ \partial_{\eta_a} \varphi = \partial_{\sigma_a} \varphi = 0 \}$. In particular, we deduce the existence of smooth functions such that \eqref{rel_const} holds near $Z$, without any singularity. 
\item If $Z$ is a space-time resonant point such that $\rho_a^Z \neq 0$, $\epsilon_1^a = \epsilon_2^a = 1$, then $\rho_b^Z = - \rho_a^Z$ and $\epsilon_1^b = \epsilon_2^b = 1$. In this case, $\xi_a+\xi_b$ vanishes at $Z$, but not $\partial_{\xi_a} \varphi$. On the sphere, the resonant point $Z$ is moreover isolated. We compute : 
\begin{align*} 
&\nabla \varphi = 24 \left( \rho_a^Z \right)^2 \begin{pmatrix} 1 \\ 1 \\ 1 \\ 1 \\ 1 \\ 1 \end{pmatrix}, \quad \nabla \partial_{\eta_a} \varphi = 6 \rho_a^Z \begin{pmatrix} -1 \\ 0 \\ 0 \\ 0 \\ 1 \\ 0 \end{pmatrix}, \quad \nabla \partial_{\sigma_a} \varphi = 6 \rho_a^Z \begin{pmatrix} 0 \\ 0 \\ -1 \\ 0 \\ 1 \\ 0 \end{pmatrix}, \\
&\quad \nabla \partial_{\eta_b} \varphi = 6 \rho_a^Z \begin{pmatrix} 0 \\ 1 \\ 0 \\ 0 \\ 0 \\ -1 \end{pmatrix}, \quad \nabla \partial_{\sigma_b} \varphi = 6 \rho_a^Z \begin{pmatrix} 0 \\ 0 \\ 0 \\ 1 \\ 0 \\ -1 \end{pmatrix} 
\end{align*}
that are linearly independent, hence the existence of a relation \eqref{rel_const} near $Z$, without singularities. 
\item If $Z$ is a space-time resonant point such that $\rho_a^Z = 0$, then $\rho_b^Z \neq 0$ (since $Z$ lies on the sphere) and $\epsilon_1^b + \epsilon_2^b + 1 \in \{ -1, 1 \}$. 

As in the previous cases, we might write in the neighborhood of $Z$ that 
\begin{align*} 
3(\eta_a + \sigma_a) (\eta_a + \rho_a) (\sigma_a + \rho_a) &= \xi_a^3 - \eta_a^3 - \sigma_a^3 - \rho_a^3 \\
&= \varphi(\xi, \eta, \sigma) + \omega \cdot \chi_1(\xi, \eta, \sigma) \partial_{\eta_b} \varphi(\xi, \eta, \sigma) \\
&\quad \quad + \omega \cdot \chi_2(\xi, \eta, \sigma) \partial_{\sigma_b} \varphi(\xi, \eta, \sigma) 
\end{align*}
by arguments of differential geometry and local inversion. 

We now try to express $\partial_{\xi_a} \varphi(\xi, \eta, \sigma)$, that can be written as : 
\[ \partial_{\xi_a} \varphi(\xi, \eta, \sigma) = (\eta_a + \sigma_a) (\eta_a + \sigma_a + 2 \rho_a) \]
in terms of
\[ (\rho_a - \eta_a) (\rho_a + \eta_a), \quad (\rho_a - \sigma_a) (\rho_a + \sigma_a), \quad (\rho_a + \eta_a)(\rho_a + \sigma_a) (\eta_a + \sigma_a) \] 
Let us set $x = \eta_a + \sigma_a$, $y = \eta_a + \rho_a$, $z = \sigma_a + \rho_a$. We search for functions $\chi_i$ such that : 
\[ x (y + z) = \chi_1 y (x - z) + \chi_2 z (x - y) + \chi_3 x y z \]
The following functions are singular solutions : $(1, 1, 2/x)$, $(1, -1, 2/y)$, $(-1, 1, 2/z)$. By applying an angular partition of unity in $(\eta_a, \sigma_a, \rho_a)$, we may combine these three solutions so to restrict the singularity to $\{ \eta_a = \sigma_a = \rho_a = 0 \}$, so of the form $\frac{1}{|\omega_a|}$. This leads to \eqref{rel_const}. 
\end{itemize} 

To conclude, we chose a partition of unity of the sphere to get \eqref{rel_const} everywhere. 
\end{Dem}

\subsection{Generic estimates} \label{section-genericweightedest}

\begin{Lem} Let $\mu$ be any symbol satisfying $\Vert \mu \Vert_{MPTT} < \infty$, and consider the following quantities : 
\begin{subequations} \label{qtlemgenestpoids}
\begin{align}
&\int_1^t \int \int e^{i s \varphi} \mu(\xi, \eta, \sigma) \widehat{f}(s, \eta) \widehat{f}(s, \sigma) \widehat{f}(s, \rho) ~ d\eta d\sigma ds \label{qtlemgenestpoids1} \\
&\int_1^t \int \int e^{i s \varphi} \mu(\xi, \eta, \sigma) \eta \widehat{f}(s, \eta) \widehat{f}(s, \sigma) \nabla_{\xi} \widehat{f}(s, \rho) ~ d\eta d\sigma ds \label{qtlemgenestpoids2} \\
&\int_1^t \int \int e^{i s \varphi} \mu(\xi, \eta, \sigma) s \sigma \cdot \nabla_{\eta}\varphi(\xi, \eta, \sigma) \widehat{f}(s, \eta) \widehat{f}(s, \sigma) \widehat{f}(s, \rho) ~ d\eta d\sigma ds \label{qtlemgenestpoids3} \\
&\int_1^t \int \int e^{i s \varphi} \mu(\xi, \eta, \sigma) s \varphi(\xi, \eta, \sigma) \widehat{f}(s, \eta) \widehat{f}(s, \sigma) \widehat{f}(s, \rho) ~ d\eta d\sigma ds \label{qtlemgenestpoids4}
\end{align}
\end{subequations}
Then for any $\alpha$, 
\[ \Vert (\ref{qtlemgenestpoids}.\alpha) \Vert_{L^2} \lesssim \Vert u \Vert_X^3 \left( 1 + \Vert u \Vert_X^2 \right) \] \label{lemgenestpoids}
\end{Lem}

\begin{Dem}
We apply theorem \ref{theobiparaproduct} and lemma \ref{lemLinfest} to obtain that, if $\delta > 0$ is small enough, 
\begin{align*}
\Vert \eqref{qtlemgenestpoids1} \Vert_{L^2} &\lesssim \int_1^t \Vert u(s) \Vert_{L^2} \Vert u(s) \Vert_{L^{\infty}}^2 ~ ds \\
&\lesssim \int_1^t \Vert u \Vert_X^3 s^{-\frac{4}{3}+2\delta} ~ ds \lesssim \Vert u \Vert_X^3 \\
\Vert \eqref{qtlemgenestpoids2} \Vert_{L^2} &\lesssim \int_1^t \Vert \nabla u(s) \Vert_{L^{\infty}} \Vert u(s) \Vert_{L^{\infty}} \Vert (x_a, x_b) f(s) \Vert_{L^2} ~ ds \\
&\lesssim \int_1^t \Vert u \Vert_X^3 s^{-\frac{11}{9}+2\delta} ~ ds \lesssim \Vert u \Vert_X^3 \\
\end{align*}
To control \eqref{qtlemgenestpoids3}, we apply an integration by parts in $\sigma$ using $e^{i s \varphi} i s \nabla_{\eta} \varphi = \nabla_{\eta} \left( e^{i s \varphi} \right)$ and obtain \eqref{qtlemgenestpoids1}, \eqref{qtlemgenestpoids2} or similar terms up to a change of variables between $\eta, \sigma, \rho$ (note that $\sigma \cdot \nabla_{\eta} \mu$ still satisfies $\Vert \sigma \cdot \nabla_{\eta} \mu \Vert_{MPTT} < \infty$). 

Finally, for \eqref{qtlemgenestpoids4}, we apply an integration by parts in time and get : 
\begin{subequations}
\begin{align}
\eqref{qtlemgenestpoids4} &= i \int_1^t \int \int e^{i s \varphi} \mu(\xi, \eta, \sigma) s \partial_s \left( \widehat{f}(s, \eta) \widehat{f}(s, \sigma) \widehat{f}(s, \rho) \right) ~ d\eta d\sigma ds \label{qtlemgenestpoids4-1} \\
&- i \int \int e^{i t \varphi} \mu(\xi, \eta, \sigma) t \widehat{f}(t, \eta) \widehat{f}(t, \sigma) \widehat{f}(t, \rho) ~ d\eta d\sigma \label{qtlemgenestpoids4-2} \\
&+ i \int \int e^{i \varphi} \mu(\xi, \eta, \sigma) \widehat{f}(1, \eta) \widehat{f}(1, \sigma) \widehat{f}(1, \rho) ~ d\eta d\sigma ds \label{qtlemgenestpoids4-3} 
\end{align}
\end{subequations}
Using lemmas \ref{lemLinfest} and \ref{lemcontrdtf}, 
\begin{align*}
\Vert \eqref{qtlemgenestpoids4-1} \Vert_{L^2} &\lesssim \int_1^t s \Vert \partial_s f(s) \Vert_{L^2} \Vert u(t) \Vert_{L^{\infty}}^2 ~ ds \\
&\lesssim \int_1^t s^{-\frac{4}{3}-\frac{1}{9}+2\delta} \Vert u \Vert_X^5 ~ ds \lesssim \Vert u \Vert_X^5 \\
\Vert \eqref{qtlemgenestpoids4-2} \Vert_{L^2} &\lesssim t \Vert u(t) \Vert_{L^2} \Vert u(t) \Vert_{L^{\infty}}^2 \\
&\lesssim t^{-\frac{1}{3}+2\delta} \Vert u \Vert_X^3 \lesssim \Vert u \Vert_X^3 
\end{align*}
And \eqref{qtlemgenestpoids4-3} is similar to \eqref{qtlemgenestpoids4-2}, replacing $t$ by $1$. 
\end{Dem}

\subsection{Proof of the weighted estimate} \label{section-proofweighted}

We develop \eqref{termtotpoids} : 
\begin{subequations}
\begin{align}
\eqref{termtotpoids} &= \partial_{\xi_a} \int_1^t \int \int e^{i s \varphi(\xi, \eta, \sigma)} (\xi_a + \xi_b) \widehat{f}(s, \eta) \widehat{f}(s, \sigma) \widehat{f}(s, \rho) ~ d\eta d\sigma \notag \\
&= \int_1^t \int \int e^{i s \varphi(\xi, \eta, \sigma)} i s (\xi_a + \xi_b) \partial_{\xi_a} \varphi(\xi, \eta, \sigma) \widehat{f}(s, \eta) \widehat{f}(s, \sigma) \widehat{f}(s, \rho) ~ d\eta d\sigma \label{termtotpoids-princ} \\
&+ \int_1^t \int \int e^{i s \varphi(\xi, \eta, \sigma)} \widehat{f}(s, \eta) \widehat{f}(s, \sigma) \widehat{f}(s, \rho) ~ d\eta d\sigma \label{termtotpoids-1} \\
&+ \int_1^t \int \int e^{i s \varphi(\xi, \eta, \sigma)} (\xi_a + \xi_b) \widehat{f}(s, \eta) \widehat{f}(s, \sigma) \partial_{\xi_a} \widehat{f}(s, \rho) ~ d\eta d\sigma \label{termtotpoids-2} 
\end{align}
\end{subequations}
Note that \eqref{termtotpoids-1} and \eqref{termtotpoids-2} are controlled by lemma \ref{lemgenestpoids}. 

To treat \eqref{termtotpoids-princ}, we apply lemma \ref{lem_rel_const} and a partition of unity to distinguish cases. 

\paragraph{Away from singularities} Assume that the relation \eqref{rel_const} holds without singularities. Let us first restrict our attention to a frequency angular area where no frequency is too large with respect to the others, that is 
\[ |\eta| \lesssim |\rho| + |\sigma|, \quad |\sigma| \lesssim |\eta| + |\rho|, \quad |\rho| \lesssim |\eta| + |\sigma| \]
Then we may apply \eqref{rel_const} and the previous inequalities to be able to control \eqref{termtotpoids-princ} by lemma \ref{lemgenestpoids} (up to some symmetry between $\eta, \sigma, \rho$). 

Consider now the case where one of the frequencies is very large with respect to the others. By symmetry, let us assume $|\rho| \gg |\eta| + |\sigma|$. Then note that $\nabla_{\eta} \varphi$ and $\nabla_{\sigma} \varphi$ do not vanish near such a point, so that 
\begin{align*}
(\xi_a + \xi_b) \partial_{\xi_a} \varphi(\xi, \eta, \sigma) &= 3 (\xi_a + \xi_b) (\xi_a - \rho_a) (\xi_a + \rho_a) \\
&= 3 (\eta_a + \sigma_a) (\xi_a + \xi_b) (\xi_a + \rho_a) \\
&= \eta_a m_1(\xi, \eta, \sigma) \cdot \nabla_{\sigma} \varphi(\xi, \eta, \sigma) + \sigma_a m_2(\xi, \eta, \sigma) \cdot \nabla_{\eta} \varphi(\xi, \eta, \sigma) 
\end{align*}
which is a relation that allows to apply lemma \ref{lemgenestpoids}. 

\paragraph{Near singularities} We now restrict our attention to a neighborhood of a point $\omega^Z = (\xi^Z, \eta^Z, \sigma^Z)$ such that $\xi_a^Z = \eta_a^Z = \sigma_a^Z = \rho_a^Z = 0$ and $\xi_b^Z = \rho_b^Z = - \eta_b^Z = - \sigma_b^Z$. Note that such singular points are isolated one from another on the sphere $\{ (\xi, \eta, \sigma), ~ |\xi|^2 + |\eta|^2 + |\sigma|^2 = 1 \}$. Let $\mu$ be a Coifman-Meyer symbol localizing near $\omega^Z$. 

Let $\chi_0 \in C^{\infty}_c(\mathbb{R}, [0, 1])$ be compactly supported in $[-2, 2]$ and identically equal to $1$ on $[-1, 1]$. Define 
\[ \chi(t, \xi, \eta, \sigma) = \chi_0\left( t^{\frac{1}{4}} \frac{|\omega_a|}{|\omega_b|} \right) \]
Note that $\Vert \chi(t, \cdot) \Vert_{MPTT} = \Vert \chi\left( 1, \cdot \right) \Vert_{MPTT} < \infty$. Let $\widetilde{\chi} := 1 - \chi$. We then separate : 
\begin{subequations}
\begin{align}
&\int_1^t \int \int e^{i s \varphi(\xi, \eta, \sigma)} i s (\xi_a + \xi_b) \partial_{\xi_a} \varphi(\xi, \eta, \sigma) \mu(\xi, \eta, \sigma) \widehat{f}(s, \eta) \widehat{f}(s, \sigma) \widehat{f}(s, \rho) ~ d\eta d\sigma \notag \\
&= \int_1^t \int \int e^{i s \varphi(\xi, \eta, \sigma)} i s (\xi_a + \xi_b) \partial_{\xi_a} \varphi(\xi, \eta, \sigma) \mu(\xi, \eta, \sigma) \chi(s, \xi, \eta, \sigma) \widehat{f}(s, \eta) \widehat{f}(s, \sigma) \widehat{f}(s, \rho) ~ d\eta d\sigma \label{termtotpoids-princ-0} \\
&+ \int_1^t \int \int e^{i s \varphi(\xi, \eta, \sigma)} i s (\xi_a + \xi_b) \partial_{\xi_a} \varphi(\xi, \eta, \sigma) \mu(\xi, \eta, \sigma) \widetilde{\chi}(s, \xi, \eta, \sigma) \widehat{f}(s, \eta) \widehat{f}(s, \sigma) \widehat{f}(s, \rho) ~ d\eta d\sigma \label{termtotpoids-princ-1}
\end{align}
\end{subequations}

We have that 
\[ \Vert |\omega_b|^{-2} \partial_{\xi_a} \varphi \chi(s) \Vert_{MPTT} \lesssim s^{-\frac{1}{2}} \]
In particular, since on the support of $\mu$ a symbol of the form $\frac{(\xi_a + \xi_b) |\omega_b|^2}{|\eta_b|^{\frac{1}{2}} |\sigma_b|^{\frac{1}{2}} |\rho_b|^{2}}$ is of Coifman-Meyer type, we may estimate using lemma \ref{lemLinfest} : 
\begin{align*} 
\Vert \eqref{termtotpoids-princ-0} \Vert_{L^2} &\lesssim \int_1^t s^{\frac{1}{2}} \Vert u(s) \Vert_{H^2} \Vert |\partial_{x_b}|^{\frac{1}{2}} u(s) \Vert_{L^{\infty}}^2 ~ ds \\
&\lesssim \int_1^t s^{\frac{1}{2}} s^{-\frac{4}{3}-\frac{1}{3}+2\delta} \Vert u \Vert_X^3 ~ ds \\
&\lesssim \Vert u \Vert_X^3
\end{align*}
choosing $\delta < \frac{1}{12}$. 

On the other hand, to control \eqref{termtotpoids-princ-1}, we apply the singular formula \eqref{rel_const} and develop : 
\begin{subequations}
\begin{align}
\eqref{termtotpoids-princ-1} &= \int_1^t \int \int e^{i s \varphi(\xi, \eta, \sigma)} i s \omega \cdot m_1^a(\xi, \eta, \sigma) \partial_{\eta_a} \varphi(\xi, \eta, \sigma) \mu(\xi, \eta, \sigma) \widetilde{\chi}(s, \xi, \eta, \sigma) \widehat{f}(s, \eta) \widehat{f}(s, \sigma) \widehat{f}(s, \rho) ~ d\eta d\sigma \label{termtotpoids-princ-1-1} \\
&+ \int_1^t \int \int e^{i s \varphi(\xi, \eta, \sigma)} i s \omega \cdot m_2^a(\xi, \eta, \sigma) \partial_{\sigma_a} \varphi(\xi, \eta, \sigma) \mu(\xi, \eta, \sigma) \widetilde{\chi}(s, \xi, \eta, \sigma) \widehat{f}(s, \eta) \widehat{f}(s, \sigma) \widehat{f}(s, \rho) ~ d\eta d\sigma \label{termtotpoids-princ-1-2} \\
&+ \int_1^t \int \int e^{i s \varphi(\xi, \eta, \sigma)} i s |\omega_a|^{-1} \omega^2 : m_1^b(\xi, \eta, \sigma) \partial_{\eta_b} \varphi(\xi, \eta, \sigma) \mu(\xi, \eta, \sigma) \widetilde{\chi}(s, \xi, \eta, \sigma) \widehat{f}(s, \eta) \widehat{f}(s, \sigma) \widehat{f}(s, \rho) ~ d\eta d\sigma \label{termtotpoids-princ-1-3} \\
&+ \int_1^t \int \int e^{i s \varphi(\xi, \eta, \sigma)} i s |\omega_a|^{-1} \omega^2 : m_2^b(\xi, \eta, \sigma) \partial_{\sigma_b} \varphi(\xi, \eta, \sigma) \mu(\xi, \eta, \sigma) \widetilde{\chi}(s, \xi, \eta, \sigma) \widehat{f}(s, \eta) \widehat{f}(s, \sigma) \widehat{f}(s, \rho) ~ d\eta d\sigma \label{termtotpoids-princ-1-4} \\
&+ \int_1^t \int \int e^{i s \varphi(\xi, \eta, \sigma)} i s |\omega_a|^{-1} \omega \cdot m_3(\xi, \eta, \sigma) \varphi(\xi, \eta, \sigma) \mu(\xi, \eta, \sigma) \widetilde{\chi}(s, \xi, \eta, \sigma) \widehat{f}(s, \eta) \widehat{f}(s, \sigma) \widehat{f}(s, \rho) ~ d\eta d\sigma \label{termtotpoids-princ-1-5} 
\end{align}
\end{subequations}
On each of these terms, we apply an integration by parts. For \eqref{termtotpoids-princ-1-1} : 
\begin{subequations}
\begin{align}
\eqref{termtotpoids-princ-1-1} &= - \int_1^t \int \int e^{i s \varphi(\xi, \eta, \sigma)} \partial_{\eta_a} \left( \omega \cdot m_1^a(\xi, \eta, \sigma) \mu(\xi, \eta, \sigma) \widetilde{\chi}(s, \xi, \eta, \sigma) \widehat{f}(s, \eta) \widehat{f}(s, \sigma) \widehat{f}(s, \rho) \right) ~ d\eta d\sigma \notag \\
&= \int_1^t \int \int e^{i s \varphi(\xi, \eta, \sigma)} |\omega_a|^{-1} \omega \cdot \widetilde{m_1^a}(\xi, \eta, \sigma) \widetilde{\mu}(\xi, \eta, \sigma) \widetilde{\chi}(s, \xi, \eta, \sigma) \widehat{f}(s, \eta) \widehat{f}(s, \sigma) \widehat{f}(s, \rho) ~ d\eta d\sigma \label{termtotpoids-princ-1-1-1} \\
&+ \int_1^t \int \int e^{i s \varphi(\xi, \eta, \sigma)} \omega \cdot \widetilde{m_1^a}(\xi, \eta, \sigma) \mu(\xi, \eta, \sigma) s^{\frac{1}{4}} |\omega_b|^{-1} \widetilde{\chi}'(s, \xi, \eta, \sigma) \widehat{f}(s, \eta) \widehat{f}(s, \sigma) \widehat{f}(s, \rho) ~ d\eta d\sigma \label{termtotpoids-princ-1-1-2} \\
&- \int_1^t \int \int e^{i s \varphi(\xi, \eta, \sigma)} \omega \cdot m_1^a(\xi, \eta, \sigma) \mu(\xi, \eta, \sigma) \widetilde{\chi}(s, \xi, \eta, \sigma) \partial_{\eta_a} \widehat{f}(s, \eta) \widehat{f}(s, \sigma) \widehat{f}(s, \rho) ~ d\eta d\sigma \label{termtotpoids-princ-1-1-3} \\
&- \int_1^t \int \int e^{i s \varphi(\xi, \eta, \sigma)} \omega \cdot m_1^a(\xi, \eta, \sigma) \mu(\xi, \eta, \sigma) \widetilde{\chi}(s, \xi, \eta, \sigma) \widehat{f}(s, \eta) \widehat{f}(s, \sigma) \partial_{\eta_a} \widehat{f}(s, \rho) ~ d\eta d\sigma \label{termtotpoids-princ-1-1-4}
\end{align}
\end{subequations}
where $\widetilde{m_1^a}$ (respectively $\widetilde{\mu}$) are symbols that are different from $m_1^a$ (respectively $\mu$) and that may vary at each line, but that does satisfy similar bounds and localisation properties. Note now that 
\[ \Vert |\omega_a|^{-1} \omega \widetilde{m_1^a} \widetilde{\mu} \widetilde{\chi}(s) \Vert_{MPTT} \lesssim s^{\frac{1}{4}} \]
where the constant is uniform in $s$. Likewise, 
\[ \Vert \omega \widetilde{m_1^a} \mu |\omega_b|^{-1} \widetilde{\chi}'(s) \Vert_{MPTT} \lesssim 1 \]
Moreover, on the support of $\mu$ (or $\widetilde{\mu}$), up to replacing $\mu$ by another Coifman-Meyer type symbol, we may replace $\omega$ by $|\sigma_b|^{\frac{1}{2}} |\rho_b|^{\frac{1}{2}}$. 
We may now estimate by lemma \ref{lemLinfest}, choosing $\delta > 0$ small enough : 
\begin{align*}
\Vert \eqref{termtotpoids-princ-1-1-1} \Vert_{L^2} &\lesssim \int_1^t s^{\frac{1}{4}} \Vert u(s) \Vert_{L^{\infty}}^2 \Vert u(s) \Vert_{L^2} ~ ds \\
&\lesssim \int_1^t s^{\frac{1}{4}} s^{-\frac{4}{3}+2\delta} \Vert u \Vert_X^3 ~ ds \lesssim \Vert u \Vert_X^3 \\
\Vert \eqref{termtotpoids-princ-1-1-3} \Vert_{L^2} &\lesssim \int_1^t \Vert |\partial_{x_b}|^{\frac{1}{2}} u(s) \Vert_{L^{\infty}}^2 \Vert x_a f \Vert_{L^2} ~ ds \\
&\lesssim \int_1^t s^{-\frac{5}{3}+2\delta} \Vert u \Vert_X^3 ~ ds \lesssim \Vert u \Vert_X^3
\end{align*}
\eqref{termtotpoids-princ-1-1-2} is similar to \eqref{termtotpoids-princ-1-1-1}, and \eqref{termtotpoids-princ-1-1-4} is symmetric to \eqref{termtotpoids-princ-1-1-3}. 

We may control \eqref{termtotpoids-princ-1-2} the same way as \eqref{termtotpoids-princ-1-1}. 

For \eqref{termtotpoids-princ-1-3}, we apply an integration by parts in $\eta_b$ : 
\begin{subequations}
\begin{align}
\eqref{termtotpoids-princ-1-3} &= - \int_1^t \int \int e^{i s \varphi(\xi, \eta, \sigma)} \partial_{\eta_b} \left( |\omega_a|^{-1} \omega^2 : m_1^b(\xi, \eta, \sigma) \mu(\xi, \eta, \sigma) \widetilde{\chi}(s, \xi, \eta, \sigma) \widehat{f}(s, \eta) \widehat{f}(s, \sigma) \widehat{f}(s, \rho) \right) ~ d\eta d\sigma \notag \\
&= \int_1^t \int \int e^{i s \varphi(\xi, \eta, \sigma)} |\omega_a|^{-1} \omega \cdot \widetilde{m_1^b}(\xi, \eta, \sigma) \widetilde{\mu}(\xi, \eta, \sigma) \widetilde{\chi}(s, \xi, \eta, \sigma) \widehat{f}(s, \eta) \widehat{f}(s, \sigma) \widehat{f}(s, \rho) ~ d\eta d\sigma \label{termtotpoids-princ-1-3-1} \\
&+ \int_1^t \int \int e^{i s \varphi(\xi, \eta, \sigma)} |\omega_a|^{-1} \omega^2 : \widetilde{m_1^b}(\xi, \eta, \sigma) \widetilde{\mu}(\xi, \eta, \sigma) s^{\frac{1}{4}} \frac{\omega_a}{|\omega_b|^2} \widetilde{\chi}'(s, \xi, \eta, \sigma) \widehat{f}(s, \eta) \widehat{f}(s, \sigma) \widehat{f}(s, \rho) ~ d\eta d\sigma \label{termtotpoids-princ-1-3-2} \\
&- \int_1^t \int \int e^{i s \varphi(\xi, \eta, \sigma)} |\omega_a|^{-1} \omega^2 : m_1^b(\xi, \eta, \sigma) \mu(\xi, \eta, \sigma) \widetilde{\chi}(s, \xi, \eta, \sigma) \partial_{\eta_b} \widehat{f}(s, \eta) \widehat{f}(s, \sigma) \widehat{f}(s, \rho) ~ d\eta d\sigma \label{termtotpoids-princ-1-3-3} \\
&- \int_1^t \int \int e^{i s \varphi(\xi, \eta, \sigma)} |\omega_a|^{-1} \omega^2 : m_1^b(\xi, \eta, \sigma) \mu(\xi, \eta, \sigma) \widetilde{\chi}(s, \xi, \eta, \sigma) \widehat{f}(s, \eta) \widehat{f}(s, \sigma) \partial_{\eta_b} \widehat{f}(s, \rho) ~ d\eta d\sigma \label{termtotpoids-princ-1-3-4}
\end{align}
\end{subequations}
We may estimate \eqref{termtotpoids-princ-1-3-1} the same way as \eqref{termtotpoids-princ-1-1-1} ; moreover, since 
\[ \Vert \omega^2 \cdot \widetilde{m_1^b}(\xi, \eta, \sigma) \widetilde{\mu}(\xi, \eta, \sigma) \frac{1}{|\omega_b|^2} \widetilde{\chi}'(s, \xi, \eta, \sigma) \Vert_{MPTT} \lesssim 1 \]
the estimate of \eqref{termtotpoids-princ-1-3-2} is similar to the estimate of \eqref{termtotpoids-princ-1-1-2}. Finally, replacing up to changing the Coifman-Meyer symbol $\omega$ by $|\sigma_b|^{\frac{1}{2}} |\rho_b|^{\frac{1}{2}}$ and using the fact that 
\[ \Vert |\omega_a|^{-1} \omega m_1^b(\xi, \eta, \sigma) \mu(\xi, \eta, \sigma) \widetilde{\chi}(s, \xi, \eta, \sigma) \Vert_{MPTT} \lesssim s^{\frac{1}{4}} \]
we may estimate 
\begin{align*}
\Vert \eqref{termtotpoids-princ-1-3-3} \Vert_{L^2} &\lesssim \int_1^t s^{\frac{1}{4}} \Vert x_b f \Vert_{L^2} \Vert |\partial_{x_b}|^{\frac{1}{2}} u(s) \Vert_{L^{\infty}}^2 ~ ds \\
&\lesssim \int_1^t s^{\frac{1}{4} - \frac{5}{3}+2\delta} \Vert u \Vert_X^3 ~ ds \lesssim \Vert u \Vert_X^3
\end{align*}
up to choosing $\delta$ small enough. The estimate of \eqref{termtotpoids-princ-1-3-4} is symmetric. 

We may control \eqref{termtotpoids-princ-1-4} the same way as \eqref{termtotpoids-princ-1-3}. 

Finally, for \eqref{termtotpoids-princ-1-5}, we apply an integration by parts in time : 
\begin{subequations}
\begin{align}
\eqref{termtotpoids-princ-1-5} &= - \int_1^t \int \int e^{i s \varphi(\xi, \eta, \sigma)} \partial_s \left( s |\omega_a|^{-1} \omega \cdot m_3(\xi, \eta, \sigma) \mu(\xi, \eta, \sigma) \widetilde{\chi}(s, \xi, \eta, \sigma) \widehat{f}(s, \eta) \widehat{f}(s, \sigma) \widehat{f}(s, \rho) \right) ~ d\eta d\sigma ds \notag \\
&+ \int \int e^{i t \varphi(\xi, \eta, \sigma)} t |\omega_a|^{-1} \omega \cdot m_3(\xi, \eta, \sigma) \mu(\xi, \eta, \sigma) \widetilde{\chi}(t, \xi, \eta, \sigma) \widehat{f}(t, \eta) \widehat{f}(t, \sigma) \widehat{f}(t, \rho) ~ d\eta d\sigma \notag \\
&- \int \int e^{i \varphi(\xi, \eta, \sigma)} |\omega_a|^{-1} \omega \cdot m_3(\xi, \eta, \sigma) \mu(\xi, \eta, \sigma) \widetilde{\chi}(1, \xi, \eta, \sigma) \widehat{f}(1, \eta) \widehat{f}(1, \sigma) \widehat{f}(1, \rho) ~ d\eta d\sigma \notag \\
&= - \int_1^t \int \int e^{i s \varphi(\xi, \eta, \sigma)} |\omega_a|^{-1} \omega \cdot m_3(\xi, \eta, \sigma) \mu(\xi, \eta, \sigma) \widetilde{\chi}(s, \xi, \eta, \sigma) \widehat{f}(s, \eta) \widehat{f}(s, \sigma) \widehat{f}(s, \rho) ~ d\eta d\sigma ds \label{termtotpoids-princ-1-5-1} \\
&- \int_1^t \int \int e^{i s \varphi(\xi, \eta, \sigma)} |\omega_a|^{-1} \omega \cdot m_3(\xi, \eta, \sigma) \mu(\xi, \eta, \sigma) s^{\frac{1}{4}} \frac{|\omega_a|}{|\omega_b|} \widetilde{\chi}'(s, \xi, \eta, \sigma) \widehat{f}(s, \eta) \widehat{f}(s, \sigma) \widehat{f}(s, \rho) ~ d\eta d\sigma ds \label{termtotpoids-princ-1-5-2} \\
&- \int_1^t \int \int e^{i s \varphi(\xi, \eta, \sigma)} s |\omega_a|^{-1} \omega \cdot m_3(\xi, \eta, \sigma) \mu(\xi, \eta, \sigma) \widetilde{\chi}(s, \xi, \eta, \sigma) \partial_s \left( \widehat{f}(s, \eta) \widehat{f}(s, \sigma) \widehat{f}(s, \rho) \right) ~ d\eta d\sigma ds \label{termtotpoids-princ-1-5-3} \\
&+ \int \int e^{i t \varphi(\xi, \eta, \sigma)} t |\omega_a|^{-1} \omega \cdot m_3(\xi, \eta, \sigma) \mu(\xi, \eta, \sigma) \widetilde{\chi}(t, \xi, \eta, \sigma) \widehat{f}(t, \eta) \widehat{f}(t, \sigma) \widehat{f}(t, \rho) ~ d\eta d\sigma \label{termtotpoids-princ-1-5-4} \\
&- \int \int e^{i \varphi(\xi, \eta, \sigma)} |\omega_a|^{-1} \omega \cdot m_3(\xi, \eta, \sigma) \mu(\xi, \eta, \sigma) \widetilde{\chi}(1, \xi, \eta, \sigma) \widehat{f}(1, \eta) \widehat{f}(1, \sigma) \widehat{f}(1, \rho) ~ d\eta d\sigma \label{termtotpoids-princ-1-5-5}
\end{align}
\end{subequations}
Note that 
\[ \Vert |\omega_a|^{-1} \omega \cdot m_3(\xi, \eta, \sigma) \mu(\xi, \eta, \sigma) \widetilde{\chi}(s, \xi, \eta, \sigma) \Vert_{MPTT} \lesssim s^{\frac{1}{4}} \]
and 
\[ \Vert |\omega_a|^{-1} \omega \cdot m_3(\xi, \eta, \sigma) \mu(\xi, \eta, \sigma) s^{\frac{1}{4}} \frac{|\omega_a|}{|\omega_b|} \widetilde{\chi}'(s, \xi, \eta, \sigma) \Vert_{MPTT} \lesssim s^{\frac{1}{4}} \]
Therefore, by lemmas \ref{lemLinfest} and \ref{lemcontrdtf}, 
\begin{align*}
\Vert \eqref{termtotpoids-princ-1-5-1} \Vert_{L^2} &\lesssim \int_1^t s^{\frac{1}{4}} \Vert u(s) \Vert_{L^2} \Vert u(s) \Vert_{L^{\infty}}^2 ~ ds \\
&\lesssim \int_1^t s^{\frac{1}{4}-\frac{4}{3}+2\delta} \Vert u \Vert_X^3 ~ ds \lesssim \Vert u \Vert_X^3 \\
\Vert \eqref{termtotpoids-princ-1-5-3} \Vert_{L^2} &\lesssim \int_1^t s^{\frac{5}{4}} \Vert \partial_s u(s) \Vert_{L^2} \Vert u(s) \Vert_{L^{\infty}}^2 ~ ds \\
&\lesssim \int_1^t s^{\frac{1}{4}-\frac{4}{3}-\frac{1}{9}+2\delta} \Vert u \Vert_X^5 ~ ds \lesssim \Vert u \Vert_X^5 \\
\Vert \eqref{termtotpoids-princ-1-5-4} \Vert_{L^2} &\lesssim t \Vert u(t) \Vert_{L^2} \Vert u(t) \Vert_{L^{\infty}}^2 \\
&\lesssim t^{-\frac{1}{3}+2\delta} \Vert u \Vert_X^3 
\end{align*}
\eqref{termtotpoids-princ-1-5-2} can be estimated the same way as \eqref{termtotpoids-princ-1-5-1}, and \eqref{termtotpoids-princ-1-5-5} the same way as \eqref{termtotpoids-princ-1-5-4} (replacing $t$ by $1$). 

This concludes the proof of proposition \ref{propaprioripoids}.


\begin{thebibliography}{0}
\bibitem{DelortAlazard_gravitywaterwaves2015} T. Alazard and J.-M. Delort, Global solutions and asymptotic behavior for two dimensional gravity water waves, {\it Ann. Sci. {\'E}c. Norm. Sup{\'e}r. (4)}, 48, No. 5 (2015), pp.~1149--1238. 

\bibitem{AnjolrasABI} P. Anjolras, Stability of the constant states in the augmented Born-Infeld system, {J. Hyperbolic Differ. Equ.}, 21, No. 4 (2024), pp.~845--948. 

\bibitem{BFR-GWPmZK2D} D. Bhattacharya, L. G. Farah, S. Roudenko, Global well-posedness for low regularity data in the 2d modified Zakharov-Kuznetsov equation, {\it J. Differential Equations}, 268, No. 12 (2020), pp.~7962--7997.

\bibitem{BiagioniLinaresLWP} H. A. Biagioni and F. Linares, Well-posedness results for the modified Zakharov-Kuznetsov equation, {\it Nonlinear equations: Methods, Models and Applications. Progress in Nonlinear Differential Equations and their applications}, 54 (2003), Birkhäuser, Basel, pp.~181-189. 

\bibitem{Masmoudietalsoliton} F. Bozgan and T.-E. Ghoul and N. Masmoudi and K. Yang, Blow-Up Dynamics for the $L^2$ critical case of the 2D Zakharov-Kuznetsov equation, {\it arXiv:2406.06568} [math.AP]. 

\bibitem{Chenetalsoliton} G. Chen and Y. Lan and X. Yuan, On the near soliton dynamics for the 2D cubic Zakharov-Kuznetsov equations, {\it arXiv:2407.00300} [math.AP]. 

\bibitem{Chironlimder} D. Chiron, Error bounds for the (KdV)/(KP-I) and (gKdV)/(gKP-I) asymptotic regime for nonlinear Schr\"odinger type equations, Ann. Inst. H. Poincar\'e{} C Anal. Non Lin\'eaire {\bf 31} (2014), no.~6, pp.~1175--1230. 

\bibitem{CoifmanMeyer} R. R. Coifman and Y. Meyer, {\it Au del{\`a} des op{\'e}rateurs pseudo-diff{\'e}rentiels}, Ast{\'e}risque, 57 (1978), Soci{\'e}t{\'e} Math{\'e}matique de France (SMF), Paris. 

\bibitem{DengIonescuPausader_EulerMaxwell} Y. Deng and A. D. Ionescu and B. Pausader, The {Euler}-{Maxwell} system for electrons: global solutions in 2D, {\it Arch. Ration. Mech. Anal.}, 225, No. 2 (2017), pp.~771--871. 

\bibitem{FaminskiiLWP} A. V. Faminskii, The Cauchy problem for the Zakharov-Kuznetsov equation, {\it Differ. Equations}, 31, No. 6 (1995), pp.~1002--1012. 

\bibitem{FarahLinaresgZK} L. G. Farah and F. Linares and A. Pastor, A note on the 2D generalized Zakharov-Kuznetsov equation: local, global, and scattering results, {\it J. Differ. Equations}, 253, No. 8 (2012), pp.~2558--2571. 

\bibitem{Germain_EulerMaxwell} P. Germain and N. Masmoudi, Global existence for the {Euler}-{Maxwell} system, {\it Ann. Sci. {\'E}c. Norm. Sup{\'e}r. (4)}, 47, No. 3 (2014), pp.~469--503. 

\bibitem{Germain_3DSchrod} P. Germain and N. Masmoudi and J. Shatah, Global solutions for 3D quadratic {Schr{\"o}dinger} equations, {\it Int. Math. Res. Not.}, 3 (2009), pp.~414--432. 

\bibitem{Germain_2DSchrod} P. Germain and N. Masmoudi and J. Shatah, Global solutions for 2D quadratic {Schr{\"o}dinger} equations, {\it J. Math. Pures Appl. (9)}, 97, No. 5 (2012), pp.~505--543. 

\bibitem{GMSWaterwaves} P. Germain and N. Masmoudi and J. Shatah, Global solutions for the gravity water waves equation in dimension 3, {\it Ann. Math. (2)}, 175 (2012), pp.~691--754. 

\bibitem{GMS_waterwavescap} P. Germain and N. Masmoudi and J. Shatah, Global existence for capillary water waves, {\it Commun. Pure Appl. Math.}, 68, No. 4 (2015), pp.~625--687. 

\bibitem{GermainPusateriRousset_mKDV} P. Germain and F. Pusateri and F. Rousset, Asymptotic stability of solitons for {mKdV}, {\it Adv. Math.}, 299 (2016), pp.~272--330. 

\bibitem{GrunrockHerrLWP} A. Grünrock and S. Herr, The Fourier restriction norm method for the Zakharov-Kuznetsov equation, {\it Discrete Contin. Dyn. Syst.}, 34, No. 5 (2014), pp.~2061--2068. 

\bibitem{GrunrockLWP3d} A. Grünrock, A remark on the modified Zakharov-Kuznetsov equation in three space dimensions, {\it Math. Res.
Lett.}, 21 (2014), pp.~127--131.

\bibitem{GrunrockgZK} A. Grünrock, On the generalized Zakharov-Kuznetsov equation at critical regularity, {\it arXiv:1509.09146}
[math.AP].

\bibitem{GuoIonescuPausader_plasma} Y. Guo and A. D. Ionescu and B. Pausader, Global solutions of certain plasma fluid models in three-dimension, {\it J. Math. Phys.}, 55, No. 12 (2014). 

\bibitem{GuoIonescuPausader_EulerMaxwell} Y. Guo and A. D. Ionescu and B. Pausader, Global solutions of the {Euler}-{Maxwell} two-fluid system in 3D, {\it Ann. Math. (2)}, 183, No. 2 (2016), pp.~377-498. 

\bibitem{GuoPausader_EulerPoisson} Y. Guo and B. Pausader, Global smooth ion dynamics in the {Euler}-{Poisson} system, {\it Commun. Math. Phys.}, 303, No. 1 (2011), pp.~89--125. 

\bibitem{NakanishiResonances} S. Gustafson and K. Nakanishi and T.-P. Tsai, Global dispersive solutions for the {Gross}-{Pitaevskii} equation in two and three dimensions, {\it Ann. Henri Poincar{\'e}}, 8, No. 7 (2007), pp.~1303--1331. 

\bibitem{NakanishiResonances2} S. Gustafson and K. Nakanishi and T.-P. Tsai, Scattering theory for the {Gross}-{Pitaevskii} equation in three dimensions, {\it Commun. Contemp. Math.}, 11, No. 4 (2009), pp.~657--707. 

\bibitem{HarropGriffLongtimemKdV} B.~H. Harrop-Griffiths, Long time behavior of solutions to the mKdV, {\it Comm. Partial Differential Equations}, 41, No.~2 (2016), pp.~282--317. 

\bibitem{HayashiNaumkinmKdV1} N. Hayashi and P.~I. Naumkin, Large time behavior of solutions for the modified Korteweg-de Vries equation, {\it Internat. Math. Res. Notices}, No.~8 (1999), pp.~395--418. 

\bibitem{HayashiNaumkinmKdV2} N. Hayashi and P.~I. Naumkin, On the modified Korteweg-de Vries equation, {\it Math. Phys. Anal. Geom.}, 4, No.~3 (2001), pp.~197--227. 

\bibitem{HerrKinoshitaLWP} S. Herr and S. Kinoshita, Subcritical well-posedness results for the Zakharov-Kuznetsov equation in dimension three and higher, {\it Ann. Inst. Fourier (Grenoble)}, 73, No. 3 (2023), pp.~1203--1267.

\bibitem{HunterIfrimTataru_ww} J. K. Hunter and M. Ifrim and D. Tataru, Two dimensional water waves in holomorphic coordinates, {\it Commun. Math. Phys.}, 346, No. 2 (2016), pp.~483--552. 

\bibitem{IfrimTataru_ww} M. Ifrim and D. Tataru, Two dimensional water waves in holomorphic coordinates. {II}: {Global} solutions, {\it Bull. Soc. Math. Fr.}, 144, No. 2 (2016), pp.~369--394. 

\bibitem{IfrimTataru1Dconj} M. Ifrim and D. Tataru, The global well-posedness conjecture for 1D cubic dispersive equations, {\it arXiv:2311.15076} [math.AP]. 

\bibitem{IonescuPausader_KG} A. D. Ionescu and B. Pausader, Global solutions of quasilinear systems of {Klein}-{Gordon} equations in 3D, {\it J. Eur. Math. Soc. (JEMS)}, 16, No. 11 (2014), pp.~2355--2431. 

\bibitem{IonescuPusateri_waterwaves} A. D. Ionescu and F. Pusateri, Global solutions for the gravity water waves system in 2d, {\it Invent. Math.}, 199, No. 3 (2015), pp.~653--804. 

\bibitem{KatoPonceCommut} T. Kato and G. Ponce, Commutator estimates and the Euler and Navier-Stokes equations, {\it Comm. Pure Appl. Math. }, 41, No.~7 (1988), pp.~891--907.

\bibitem{KinoshitaGWP} S. Kinoshita, Global well-posedness for the Cauchy problem of the Zakharov-Kuznetsov equation in 2D, {\it Ann. Inst. Henri Poincaré, Anal. Non Linéaire}, 38, No. 2 (2021), pp.~451--505. 

\bibitem{KinoshitaLWPmZK2D} S. Kinoshita, Well-posedness for the Cauchy problem of the modified Zakharov-Kuznetsov equation, {\it Funkcial. Ekvac.}, 65, No. 2 (2022), pp.~139--158. 

\bibitem{LLS-derivation} D. Lannes, F. Linares and J.-C. Saut, The Cauchy problem for the Euler-Poisson system and derivation of the Zakharov-Kuznetsov equation, {\it Studies in phase space analysis with applications to PDEs}, 181--213, Progr. Nonlinear Differential Equations Appl., 84, Birkh\"auser/Springer, New York, pp.~181--213. 

\bibitem{LinaresPastorLWP} F. Linares and A. Pastor, Well-posedness for the two-dimensional modified Zakharov-Kuznetsov equation, {\it SIAM J. Math. Anal.}, 41, No. 4 (2009), pp.~1323--1339. 

\bibitem{LinaresPastorgZK} F. Linares and A. Pastor, Local and global well-posedness for the 2D generalized Zakharov-Kuznetsov equation, {\it J. Funct. Anal.}, 260, No. 4 (2011), pp.~1060-1085. 

\bibitem{LinaresPoncedispersivebook} F. Linares and G. Ponce, {\it Introduction to nonlinear dispersive equations.} New York, NY: Springer (2009). 

\bibitem{LinaresSautLWP} F. Linares and J.-C. Saut, The Cauchy problem for the 3D Zakharov-Kuznetsov equation, {\it Discrete Contin.
Dyn. Syst.}, 24 (2009), pp.~547--565.

\bibitem{MolinetPilodLWP} L. Molinet and D. Pilod, Bilinear Strichartz estimates for the Zakharov-Kuznetsov equation and applications, {\it Ann. Inst. Henri Poincaré, Anal. Non Linéaire}, 32, No. 2 (2015), pp.~347--371. 

\bibitem{MPTTBiparameter2004} C. Muscalu and J. Pipher and T. Tao and C. Thiele, Bi-parameter paraproducts, {\it Acta Math.}, 193, Vol. 2 (2004), pp.~269-296. 

\bibitem{MuscaluHAbook1} C. Muscalu and W. Schlag, {\it Classical and multilinear harmonic analysis. Volume I.}, Cambridge University Press (2013). 

\bibitem{MuscaluHAbook2} C. Muscalu and W. Schlag, {\it Classical and multilinear harmonic analysis. Volume II.}, Cambridge University Press (2013). 

\bibitem{ShatahPusateri} F. Pusateri and J. Shatah, Space-time resonances and the null condition for first-order systems of wave equations, {\it Commun. Pure Appl. Math.}, 66, No. 10 (2013), pp.~1495--1540. 

\bibitem{RibaudVentoLWP} F. Ribaud and S. Vento, A note on the Cauchy problem for the 2D generalized Zakharov-Kuznetsov equations, {\it C. R., Math., Acad. Sci. Paris}, 350, No. 9--10 (2012), pp.~499--503. 

\bibitem{RibaudVentoLWP3d} F. Ribaud, S. Vento, Well-posedness results for the three-dimensional Zakharov-Kuznetsov equation, {\it SIAM J. Math. Anal.}, 44 (2012), pp.~2289--2304.

\bibitem{SteinFourierbook} E. M. Stein and G. Weiss, {\it Introduction to Fourier analysis on Euclidean spaces}, Princeton Math. Ser. 32 (1971), Princeton University Press, Princeton, NJ. 

\bibitem{ZKoriginal} V. E. Zakharov and E. A. Kuznetsov, Three-dimensional solitons, {\it Sov. Phys. JETP}, 39 (1974), pp.~285--286.
\end{thebibliography}
\end{document}